\documentclass{amsproc}

\usepackage{amsmath,amssymb}
\newtheorem{theo}{Theorem}[section]
\newtheorem{lemm}[theo]{Lemma}
\newtheorem{prop}[theo]{Proposition}
\newtheorem{defi}[theo]{Definition}

\newtheorem{coro}[theo]{Corollary}
\theoremstyle{remark}
\newtheorem{remark}[theo]{Remark}

\numberwithin{equation}{section}

\def\cen{\centerline}

\def\om{\omega}

\def\a{\alpha}

\def\mod{\hbox{mod}}

\newfont{\df}{eufm10}
\def\vep{\varepsilon}
\def\ep{\epsilon}

\def\ot{\otimes}

\def\ot{\otimes}

\def\ra{\rangle}
\def\la{\langle}
\def\lg{\langle}
\def\rg{\rangle}

\def\be{\beta}

\def\ep{\epsilon}
\def\ot{\otimes}

\def\om{\omega}

\begin{document}

\title[two-parameter quantum affine algebra for type $C_n^{(1)}$]
{Two-parameter Quantum Affine Algebra of Type ${\mathrm C_n^{(1)}},$
Drinfeld Realization and Vertex Representation}


\author[Hu]{Naihong Hu}
\address{Department of Mathematics,  Shanghai Key Laboratory of Pure Mathematics and Mathematical Practice, East China Normal University,
Shanghai 200241, PR China} \email{nhhu@math.ecnu.edu.cn}

\author[Zhang]{Honglian Zhang$^{\star}$}
\address{Department of Mathematics,
Shanghai University, Shanghai 200444, PR China}
\email{hlzhangmath@shu.edu.cn}

\thanks{$^\star$ H. Zhang is the corresponding author of this paper}

%    General info
%\date{Sept. 1, 2004 and, in revised form, May 18, 2005.}


%%%%%%%%%%%%%%%%%%%%%%%%%%%%%%%%%%%%%%%%%%%%%%%%%%%%%%%%%%%%%%%%%%%%%%%%
%\footnote{Corresponding author.}%
%%%%%%%%%%%%%%%%%%%%%%%%%%%%%%%%%%%%%%%%%%%%%%%%%%%%%%%%%%%%%%%%%%%%%%%%
\begin{abstract}
The two-parameter quantum vertex operator representation of level-one is explicitly constructed for $U_{r,s}(C^{(1)}_n)$ based on its two-parameter Drinfeld realization we give.
This construction will degenerate to the one-parameter case due to Jing-Koyama-Misra (\cite{JKM2}) when $rs=1$.
\end{abstract}

\maketitle

%\section*{This is an unnumbered first-level section head}
%This is an example of an unnumbered first-level heading.

%\specialsection*{This is a Special Section Head}
%This is an example of a special section head%

\section{Introduction}
\medskip

In 2000, the study of two-parameter quantum groups was revitalized by a series of work for type $A$ of Benkart and Witherspoon \cite{BW1, BW2, BW3} originally
obtained by Takeuchi \cite{T}. A systematic study afterwards both on the structures and finite-dimensional representation theory of two-parameter quantum groups for the semisimple Lie algebras of  any other types can be seen \cite{BGH1, BGH2}, \cite{BH}, \cite{HP1, HP2}, \cite{HS, HS1, HW1, HW2}, etc.
In 2004, Hu, Rosso and Zhang \cite{HRZ} began to investigate the two-parameter
quantum affine algebra of type $A_n^{(1)}$, and obtained
the two-parameter version of the celebrated Drinfeld realization in the case of $U_{r,s}(\widehat{\frak{sl}_n})$, as well as proposed for the first time the quantum affine
Lyndon basis as a monomial basis in the affine case. A general insight (\cite{HZ}) for handling the two-parameter quantum affine algebras of untwisted types in a unified manner had been found
when the first author visited ICTP early in 2006, that is, the $\tau$-invariant generating functions for the two-parameter version (where $\tau$ is the involution as a $\mathbb Q$-antiautomorphism) successfully served as a defining tool of the Drinfeld realization formalism in a compact form avoiding the case-by-case manner. As a valid verification of such defining relations for
Drinfeld realizations, the quantum two-parameter vertex representations of level one for the simply-laced cases $X^{(1)}_n$, where $X=A, D, E$, had been established there.
Also for type $G_2^{(1)}$, the validness of our definition for the two-parameter Drinfeld realization can be well-checked in the level of its two-parameter quantum vertex representation,
see \cite{GHZ}. As two generalizations of \cite{HZ} to the two-parameter quantum affine algebras of twisted types $X^{(r)}$ for $r=2,3$, the readers can consult \cite{JZ1},
to the two-parameter quantum toroidal algebras, please refer to \cite{JZ2}, where the authors gave a McKay correspondence formalism of the vertex representations obtained in \cite{HZ}.

It was known that the theory of two-parameter quantum affine algebras has been developed with some analogous stories as in the
one-parameter counterpart such as Drinfeld realization theorem with a different argument approach \cite{HRZ, HZ, JZ1}, fermionic realization \cite{JZ3} and a finite-dimensional representation theory \cite{JZ4}.
The quantum vertex representations of one-parameter quantum affine algebras for the untwisted types were
first constructed by Frenkel-Jing \cite{FJ} that confirmed the Drinfeld's celebrated conjectural ``new realization" \cite{Dr2} in the level of vertex representations, even though
it was not proved until Beck gave his rigorous proof for the untwisted types by generalizing the Lusztig's braid automorphisms suitable for the quantum affine algebras based on the work of Damiani \cite{D} and Levendorskii-Soibel\'man-Stukopin \cite{LSS}.  Afterwards, the quantum vertex representation theory had been established in many works, for instance, see
\cite{Be, J1, JM, J2, JKM1, JKM2} and references therein, etc.

The goal of the current paper is to construct the level-one vertex representation
of two-parameter quantum affine algebra of ${U}_{r,s}(\mathrm{C}_n^{(1)})$, which also verifies our unified defining formalism for the two-parameter version of the Drinfeld realization
for the multiply-laced cases in the level of two-parameter quantum vertex representations.

The paper is organized as follows. In section 2, We first give Drinfeld-Jimbo presentation of
two-parameter quantum affine algebra ${U}_{r,s}(\mathrm{C}_n^{(1)})$ in the sense of Hopf algebra.
The Drinfeld realization of two-parameter quantum affine algebra
${U}_{r,s}(\mathrm{C}_n^{(1)})$ is given in section 3.
Furthermore we present and prove the Drinfeld theorem between the above two realizations.
In section 4, we start from the two-parameter enlarged quantum Heisenberg algebra and introduce a quasi-cocycle. Then we
construct the level-one quantum vertex representation of two-parameter
quantum affine algebras ${U}_{r,s}(\mathrm{C}_n^{(1)})$. This construction will degenerate to the one-parameter case due to Jing-Koyama-Misra \cite{JKM2} when $rs=1$.

\section{Quantum Affine Algebra ${U}_{r,s}(\mathrm{C}_n^{(1)})$ and Drinfeld Double}

\subsection{Structure of $U_{r,s}(\mathrm{C}_n^{(1)})$}\,
%From now on we restrict ourselves to the case of
%$\mathrm{C}_n^{(1)}$ unless stated otherwise.
Let $\mathbb{K}=\mathbb{Q}(r^{\frac1{2}},s^{\frac1{2}})$ denote a field of rational functions
with two-parameters $r^{\frac1{2}},\,s^{\frac1{2}}$ ($r\ne \pm s$).
Let $\frak{g}$ be the finite-dimensional complex simple Lie algebra of type $C_n$,
  with Cartan subalgebra $\frak{h}$ and Cartan matrix $A=(a_{ij})_{i,j\in I}$ $(I=\{1, 2,\cdots, n\})$.
Fix coprime integers
  $(d_i)_{i\in I}$ such that $(d_ia_{ij})$ is symmetric.
Assume $\Phi$ is a finite root system of type $C_n$ with
$\Pi$ a base of simple roots. Regard $\Phi$ as a subset of a
Euclidean space E = ${\mathbb R}^n$ with an inner product $\langle
\, ,\, \rangle$. Let $\epsilon _{1},\epsilon _{2}, \cdots, \epsilon_{n}$
denote an orthonormal basis of E.
    Let $\alpha_{i}
= \epsilon_{i}-\epsilon_{i+1},\,\alpha_{n} =2\epsilon_{n}$ be the simple roots of the simple Lie algebra $\frak{sp}_{2n}$,
$\{ \alpha_i^\vee \}$ and $\{ \lambda_i\}$, the sets of, respectively, simple coroots and fundamental weights.
$Q = \bigoplus_{i=1}^{n}\mathbb{Z}\alpha_i$
is the root lattice.
Let $\theta$ be the highest root and $\delta$ denote the primitive imaginary root of $C_n^{(1)}$. Take
$\alpha_0=\delta-\theta$, then $\Pi'=\{\alpha_i\mid
i\in I_0=\{0,\,1,\,\cdots, n\}\}$ is a base of simple roots of the affine Lie
algebra ${C_n^{(1)}}$. Let
$\widehat{Q}=\bigoplus_{i=0}^{n}\mathbb{Z}\alpha_i$ denote the root lattice of ${C_n^{(1)}}$.

Let $c$ be the canonical central element of the affine Lie algebra of
type $C_n^{(1)}$. Define a  nondegenerate
symmetric bilinear form $(\,  \mid \,)$ on ${\frak{h^*}}$ satisfying
$$(\alpha_i\mid \alpha_j)=d_ia_{ij}, \qquad (\delta \mid \alpha_i)=(\delta\mid \delta)=0, \quad
\textit {for all } \ i,\,j \in I_0,$$ where
$(d_0,\,d_1,\,\cdots, d_n)=(1,\,1/2,\,\cdots, 1/2,\, 1)$. Denote $r_i = r^{d_i}$, $s_i =
s^{d_i}$ for $i=0,\,1,\,\cdots, n$.

Given two sets of symbols $W=\{\omega_0,\omega_1, \cdots,\omega_n\}$,
$W'=\{\omega_0',\omega_1', \cdots, \omega_n'\}$. Define the structural
constants matrix $(\langle \omega_i',\omega_j \rangle)_{(n+1)\times
{n+1}}$ of type $C_n^{(1)}$ by
$$\left(\begin{array}{cccccc}
rs^{-1}& r^{-1}& 1 & \cdots & 1 & rs \\
s & r^{\frac{1}{2}}s^{-\frac{1}{2}} & r^{-\frac{1}{2}}  & \cdots & 1 & 1\\
\cdots &\cdots &\cdots & \cdots & \cdots & \cdots\\
1 & 1 & 1  & \cdots & r^{\frac{1}{2}}s^{-\frac{1}{2}} & r^{-1}\\
 (rs)^{-1} & 1 & 1 & \cdots & s & rs^{-1}
\end{array}\right)$$

%For simple Lie algebra
%$\mathrm{C}_n^{(1)}$, ~Fix $\Pi=\{\al_i=\ep_i-\ep_{i+1}\mid
%1\leqslant i<n\}\cup\{\al_n=2\ep_n\}$ and
%$\Phi=\{\pm\ep_i\pm\ep_j\mid 1\leqslant i\ne j\leqslant n\}\cup\{\pm
%2\ep_i\mid 1\leqslant i\leqslant n\}$. ~The highest positive root
%$\theta=2\alpha_1+\cdots+2\alpha_{n-1}+\alpha_n$, ~so the null root
%$\alpha_0=\delta-2\epsilon_1$. ~In this case, let
%$r_i=r^{\frac{(\al_i, \al_i)}{2}}$ ºÍ $s_i=s^{\frac{(\al_i,
%\al_i)}{2}}$, ~then
% $r_0=r_n=r,\,~r_1=\cdots=r_{n-1}=r^{\frac{1}{2}}$ and $s_0=s_n=s,\,
%~s_1=\cdots=s_{n-1}=s^{\frac{1}{2}}$.

\begin{defi} {\it The {\bf two-parameter quantum affine algebra
 $U_{r,s}(\mathrm{C}_n^{(1)})$} is a unital associative
algebra over $\mathbb{K}$ generated by the elements $e_j,\, f_j,\,
\omega_j^{\pm 1},\, \omega_j'^{\,\pm 1}\, (j\in I_0),\,
\gamma^{\pm\frac{1}2},\,\gamma'^{\pm\frac{1}2},$ $D^{\pm1},
D'^{\,\pm1}$, satisfying the following relations:
\medskip

\noindent
$(\hat{C}1)$
$\gamma^{\pm\frac{1}2},\,\gamma'^{\pm\frac{1}2}$
are central with $\gamma=\om_\delta$, $\gamma'=\om'_\delta$, such
that $\gamma\gamma'=(rs)^c$. The $\omega_i^{\pm 1}, {\omega'}_j^{\,\pm
1}$ all commute with one another and
$\omega_i\,\omega_i^{-1}=\omega_i'\,\omega_i'^{\,-1}=1$,
$[\,\om_i^{\pm1}, D^{\pm1}\,]=[\,\om_j'^{\,\pm1},
D^{\pm1}\,]=[\,\om_i^{\pm1}, D'^{\pm1}\,]=[\,\om_j'^{\,\pm1},
D'^{\pm1}\,]=[D'^{\,\pm1},
D^{\pm1}]=0$.\\
%\medskip

\noindent
$(\hat{C}2)$ \ \textit{For} $\,0 \leqslant i \leqslant n$
\textit{and} $1\leqslant j<n$,
\begin{equation*}\begin{array}{lll}
&D\,e_i\,D^{-1}=r_i^{\delta_{0i}}\,e_i,\qquad\qquad\quad
&D\,f_i\,D^{-1}=r_i^{-\delta_{0i}}\,f_i,\\
&\omega_je_i\omega_j^{-1}=r_j^{(\epsilon_j,
\alpha_i)}s_j^{(\epsilon_{j+1}, \alpha_i )}\,e_i,~~~~~~~~~~~
&\omega_jf_i\omega_j^{-1}=r_j^{-(\epsilon_j,
\alpha_i)}s_j^{-(\epsilon_{j+1},\alpha_i)}\,f_i,\\
&\omega_n e_j\omega_n^{-1}=r^{(\epsilon_n,
\alpha_j)}\,e_j,
&\omega_n f_j\omega_n^{-1}=r^{-(\epsilon_n,
\alpha_j)}\,f_j,\\
&\omega_ne_n\omega_n^{-1}=r^{\frac{1}{2}(\epsilon_n,
\alpha_n)}s^{-\frac{1}{2}(\epsilon_n, \alpha_n)}\,e_n,~~~~~~~~~
&\omega_nf_n\omega_n^{-1}=r^{-\frac{1}{2}(\epsilon_n,
\alpha_n)}s^{\frac{1}{2}(\epsilon_n, \alpha_n)}\,f_n,\\
&\omega_0e_j\omega_0^{-1}=r^{-(\epsilon_{j+1}, \alpha_0)}
s^{(\epsilon_1, \alpha_j)}\,e_j,
& \omega_0f_j\omega_0^{-1}=r^{(\epsilon_{j+1}, \alpha_0 )}
s^{-(\epsilon_1, \alpha_j)}\,f_j ,\\
&\omega_n e_0\omega_n^{-1}=rs\, e_0,~~~~~~~
&\omega_n f_0\omega_n^{-1}=(rs)^{-1}\,f_0 ,\\
&\omega_0e_n\omega_0^{-1}=(rs)^{-1}\,e_n,~~~~~~
&\omega_0f_n\omega_0^{-1}=rs \,f_0 ,\\
&\omega_0e_0\omega_0^{-1}=r^{-\frac{1}{2}(\epsilon_1,
\alpha_0)}s^{\frac{1}{2}(\epsilon_1, \alpha_0)}e_0 ~~~~~~
&\omega_0f_0\omega_0^{-1}=r^{\frac{1}{2}(\epsilon_1,
\alpha_0)}s^{-\frac{1}{2}(\epsilon_1, \alpha_0)}f_0 .
\end{array}\end{equation*}
$(\hat{C}3)$ \ \textit{For} $\,0 \leqslant i \leqslant n$
\textit{and}
 $1\leqslant j<n$,
\begin{equation*}\begin{array}{lll}
 &D'\,e_i\,D'^{-1}=s_i^{\delta_{0i}}\,e_i,\qquad\qquad\quad
&D'\,f_i\,D'^{-1}=s_i^{-\delta_{0i}}\,f_i,\\
&\omega'_je_i{\omega'}_j^{-1}=s_j^{(\epsilon_j,
\alpha_i)}r_j^{(\epsilon_{j+1}, \alpha_i )}\,e_i,~~~~~~~~~~~
&\omega'_jf_i{\omega'}_j^{-1}=s_j^{-(\epsilon_j,
\alpha_i)}r_j^{-(\epsilon_{j+1},\alpha_i)}\,f_i,\\
&\omega'_n e_j{\omega'}_n^{-1}=s^{(\epsilon_n,
\alpha_j)}\,e_j,
&\omega'_n f_j{\omega'}_n^{-1}=s^{-(\epsilon_n,
\alpha_j)}\,f_j,\\
&\omega'_n e_n{\omega'}_n^{-1}=s^{\frac{1}{2}(\epsilon_n,
\alpha_n)}r^{-\frac{1}{2}(\epsilon_n, \alpha_n)}\,e_n,~~~~~~~~~
&\omega'_nf_n{\omega'}_n^{-1}=s^{-\frac{1}{2}(\epsilon_n,
\alpha_n)}r^{\frac{1}{2}(\epsilon_n, \alpha_n)}\,f_n,\\
&\omega'_0e_j{\omega'}_0^{-1}=s^{-(\epsilon_{j+1}, \alpha_0)}
r^{(\epsilon_1, \alpha_j)}\,e_j
& \omega'_0f_j{\omega'}_0^{-1}=s^{(\epsilon_{j+1}, \alpha_0 )}
r^{-(\epsilon_1, \alpha_j)}\,f_j ,\\
&\omega'_n e_0{\omega'}_n^{-1}=rs\, e_0,~~~~~~~
&\omega'_n f_0{\omega'}_n^{-1}=(rs)^{-1}\,f_0 ,\\
&\omega'_0e_n{\omega'}_0^{-1}=(rs)^{-1}\,e_n, ~~~~~~
&\omega'_0f_n{\omega'}_0^{-1}=rs \,f_0 ,\\
&\omega'_0e_0{\omega'}_0^{-1}=s^{-\frac{1}{2}(\epsilon_1,
\alpha_0)}r^{\frac{1}{2}(\epsilon_1, \alpha_0)}e_0 ~~~~~~
&\omega'_0f_0{\omega'}_0^{\,-1}=s^{\frac{1}{2}(\epsilon_1,
\alpha_0)}r^{-\frac{1}{2}(\epsilon_1, \alpha_0)}f_0 .
\end{array}\end{equation*}
$(\hat{C}4)$ \ \textit{For} $\,i,\, j\in I_0$, \textit{~then we
have:}
 $$[\,e_i, f_j\,]=\frac{\delta_{ij}}{r_i-s_i}(\omega_i-\omega'_i).$$
$(\hat{C}5)$
  \textit{For all}  $1\leqslant i\ne j\leqslant n$ \textit{but}
  $(i, j)\not\in\{(0, n), (n, 0)\}$ \textit{such that} $a_{ij}=0$,  \textit{~then we have:}
 $$[e_i, e_j]=[f_i, f_j]=0,$$
 $$e_ne_0=rs\,e_0e_n,\qquad f_0f_n=rs\,f_nf_0.$$
$(\hat{C}6)$ \ \textit{For} $1\leqslant i\leqslant n-2$,
\textit{~there are the following }$(r,s)$-Serre \textit{relations:}
\begin{gather*}
e_i^2e_{i+1}-(r{+}s)\,e_ie_{i+1}e_i+(rs)\,e_{i+1}e_i^2=0,\\
e_0^2e_{1}-(r{+}s)\,e_0e_{1}e_0+rs\,e_{1}e_0^2=0,\\
e_{i+1}^2e_i-(r_{i+1}^{-1}{+}s_{i+1}^{-1})\,e_{i+1}e_ie_{i+1}+
(r_{i+1}^{-1}s_{i+1}^{-1})\,e_ie_{i+1}^2=0,\\
e_n^2e_{n-1}-(r^{-1}{+}s^{-1})\,e_ne_{n-1}e_n+(r^{-1}s^{-1})\,e_{n-1}e_n^2=0,\\
e_{n-1}^3e_n-(r{+}(rs)^{\frac{1}{2}}{+}s)\,e_{n-1}^2e_ne_{n-1}+(rs)^{\frac{1}{2}}\,
(r{+}(rs)^{\frac{1}{2}}{+}s)\,e_{n-1}e_ne_{n-1}^2-
(rs)^{\frac{3}{2}}\,e_ne_{n-1}^3=0,\\
e_{1}^3e_0-(r^{-1}{+}(rs)^{-\frac{1}{2}}{+}s^{-1})\,e_{1}^2e_0e_{1}+
(rs)^{-\frac{1}{2}}\,
(r^{-1}{+}(rs)^{-\frac{1}{2}}{+}s^{-1})\,e_{1}e_0e_{1}^2-
(rs)^{-\frac{3}{2}}\,e_0e_{1}^3=0.
\end{gather*}
 $(\hat{C}7)$ \textit{For} $1\leqslant i\leqslant n-2$, \textit{~there are the following }
$(r,s)$-Serre \textit{relations:}
\begin{gather*}
f_{i+1}f_i^2-(r{+}s)\,f_if_{i+1}f_i+(rs)\,f_i^2f_{i+1}=0,\\
f_{1}f_0^2-(r{+}s)\,f_0f_{1}f_0+rs\,f_0^2f_{1}=0,\\
f_if_{i+1}^2-(r_{i+1}^{-1}{+}s_{i+1}^{-1})\,f_{i+1}f_if_{i+1}+
(r_{i+1}^{-1}s_{i+1}^{-1})\,f_{i+1}^2f_i=0,\\
f_{n-1}f_n^2-(r^{-1}{+}s^{-1})\,f_nf_{n-1}f_n+(r^{-1}s^{-1})\,f_n^2f_{n-1}=0,\\
f_nf_{n-1}^3-(r{+}(rs)^{\frac{1}{2}}{+}s)\,f_{n-1}f_nf_{n-1}^2+(rs)^{\frac{1}{2}}\,
(r{+}(rs)^{\frac{1}{2}}{+}s)\,f_{n-1}^2f_nf_{n-1}-
(rs)^{\frac{3}{2}}\,f_{n-1}^3f_n=0,\\
f_0f_{1}^3-(r^{-1}{+}(rs)^{-\frac{1}{2}}{+}s^{-1})\,f_{1}f_0f_{1}^2+
(rs)^{-\frac{1}{2}}\,
(r^{-1}{+}(rs)^{-\frac{1}{2}}{+}s^{-1})\,f_{1}^2f_0f_{1}-
(rs)^{-\frac{3}{2}}\,f_{1}^3f_0=0.
\end{gather*}}

%\begin{remark} The relations $(\hat{C}2)$ and
% $(\hat{C}3)$ can be determined by the following two-parameter quantum Cartan matrix
% of type $\mathrm{C}_{n}^{(1)}$.
%$$\left(\begin{array}{cccccc}
%rs^{-1}& r^{-1}& 1 & \cdots & 1 & rs \\
%s & r^{\frac{1}{2}}s^{-\frac{1}{2}} & r^{-\frac{1}{2}}  & \cdots & 1 & 1\\
%\cdots &\cdots &\cdots & \cdots & \cdots & \cdots\\
%1 & 1 & 1  & \cdots & r^{\frac{1}{2}}s^{-\frac{1}{2}} & r^{-1}\\
% (rs)^{-1} & 1 & 1 & \cdots & s & rs^{-1}
%\end{array}\right)$$
%\end{remark}

%\smallskip
 $U_{r,s}(\mathrm{C}_n^{(1)})$ is a Hopf algebra with the coproduct
$\Delta$, the counit $\vep$ and the antipode $S$ defined below: for
$i\in I_0$, we have
$$\Delta(\gamma^{\pm\frac{1}2})=\gamma^{\pm\frac{1}2}\otimes
\gamma^{\pm\frac{1}2}, \qquad
\Delta(\gamma'^{\,\pm\frac{1}2})=\gamma'^{\,\pm\frac{1}2}\otimes
\gamma'^{\,\pm\frac{1}2},$$
$$\Delta(D^{\pm1})=D^{\pm1}\otimes
D^{\pm1},\qquad \Delta(D'^{\,\pm1})=D'^{\,\pm1}\otimes
D'^{\,\pm1},$$
$$\Delta(\om_i)=\om_i\ot \om_i, \qquad
\Delta(\om_i')=\om_i'\ot \om_i',$$
$$\Delta(e_i)=e_i\ot 1+\om_i\ot e_i, \qquad \Delta(f_i)=f_i\ot \om_i'+1\ot f_i, $$
$$\vep(e_i)=\vep(f_i)=0,\quad \vep(\gamma^{\pm\frac{1}2})
=\vep(\gamma'^{\,\pm\frac{1}2})=\vep(D^{\pm1})=\vep(D'^{\,\pm1})=\vep(\om_i)=\vep(\om_i')=1,
$$
$$S(\gamma^{\pm\frac{1}2})=\gamma^{\mp\frac{1}2},\quad S(\gamma'^{\pm\frac{1}2})=\gamma'^{\mp\frac{1}2},\quad
S(D^{\pm1})=D^{\mp1},\quad S(D'^{\,\pm1})=D'^{\,\mp1},$$
$$S(e_i)=-\om_i^{-1}e_i,\quad
S(f_i)=-f_i\,\om_i'^{-1},\quad S(\om_i)=\om_i^{-1}, \quad
S(\om_i')=\om_i'^{-1}.$$
\end{defi}

\subsection{Triangular decomposition of $
U_{r,s}(\mathrm{C}_n^{(1)})$}
\begin{coro}
$U_{r,s}(\mathrm{C}_n^{(1)})\cong U_{r,s}(\widehat{\frak n}^-)\ot
U^0\ot U_{r,s}(\widehat{\frak n})$, as vector spaces,
where $U^0=\mathbb
K[\om_0^{\pm1},\om_1^{\pm1},\cdots,\om_n^{\pm1},{\om_0'}^{\pm1},$
${\om_1'}^{\pm1},\cdots,{\om_n'}^{\pm1}]$, and $U_{r,s}(\widehat{\frak n})$ $($resp.
$U_{r,s}(\widehat{\frak n}^-)$\,$)$ is the subalgebra generated
by $e_i$ $($resp. $f_i$$)$ for all $i\in I_0$.
$\hat{\mathcal B}$ $($resp. $\hat{\mathcal B'}$$)$ is the Borel Hopf
subalgebra of $U_{r, s}(\mathrm{C}_n^{(1)})$ generated by $e_j,
\omega_j^{ \pm1}$, $\gamma^{\pm\frac{1}2}, D^{\pm1}$ $($resp. $f_j,\,
\omega_j'^{\, \pm1}, \, \gamma'^{\pm\frac{1}2},\, D'^{\,\pm1})$ with
$j\in I_0$.
\end{coro}
\begin{defi}
Let $\tau$ be the $\mathbb{Q}$-algebra anti-automorphism of
$U_{r,s}(C_n^{(1)})$ such that $\tau(r)=s$, $\tau(s)=r$, $\tau(\la
\om_i',\om_j\ra^{\pm1})=\la \om_j',\om_i\ra^{\mp1}$, and
\begin{gather*}
\tau(e_i)=f_i, \quad \tau(f_i)=e_i, \quad \tau(\om_i)=\om_i',\quad
\tau(\om_i')=\om_i,\\
\tau(\gamma)=\gamma',\quad
\tau(\gamma')=\gamma,\quad\tau(D)=D',\quad \tau(D')=D.
\end{gather*}
Then ${\widehat{\mathcal B'}}=\tau({\widehat{\mathcal B}})$ with
those induced defining relations from ${\widehat{\mathcal B}}$,
those cross relations in $(\hat{\rm C2})$---$(\hat{\rm C4})$,
$(\hat{\rm C5})$ and $(\hat{\rm C6})$ are antisymmetric with respect
to $\tau$.
\end{defi}

\section{Drinfeld Realization of  $U_{r,s}(C_n^{(1)})$}
\subsection{Drinfeld Realization} In this subsection, we
describe the two-parameter Drinfeld quantum affinization of
$U_{r,s}(C_n)$, which will play an important role in the vertex
representation of Section 4. Briefly denote $\la i,
j\ra:=\la\omega_i', \omega_j\ra$.

\begin{defi} Let ${\mathcal U}_{r,s}(C_n^{(1)})$  be the unital
associative algebra over $\mathbb{K}$ generated by the elements
$x_i^{\pm}(k)$, $a_i(\ell)$, $\om_i^{\pm1}$, ${\om'}_i^{\pm1}$,
$\gamma^{\pm\frac{1}{2}}$, ${\gamma'}^{\,\pm\frac{1}2}$, $D^{\pm1}$,
$D'^{\,\pm1}$ $(i=1,\,\cdots, n$, $k,\,k' \in \mathbb{Z}$, $l,\,l' \in
\mathbb{Z}\backslash \{0\})$ with the following defining relations:
\begin{eqnarray*}
({\rm D1})&&\gamma^{\pm\frac{1}{2}}, \gamma'^{\,\pm\frac{1}{2}}
\textit{ are central with } \gamma\gamma'=(rs)^{c}\,\textit{and }
\omega_i\,\omega_i^{-1}=\omega_i'\,\omega_i'^{\,-1}=1 \
(i\in I),\\
&& [\,\omega_i^{\pm 1}, \omega_j^{\,\pm 1}\,]=[\,\omega_i^{\pm 1},
\omega_j'^{\,\pm 1}\,]=[\,\omega_i'^{\pm 1}, \omega_j'^{\,\pm
1}\,]=[\,\om_i^{\pm1}, D^{\pm1}\,]=[\,\om_j'^{\,\pm1},
D^{\pm1}\,]\\
&&=[\,\om_i^{\pm1}, D'^{\pm1}\,]=[\,\om_j'^{\,\pm1},
D'^{\pm1}\,]=[D'^{\,\pm1}, D^{\pm1}]=0, \ \textit{for } \, i,\,j\in I.\\
({\rm D2})
 &&[a_i(\ell),a_j(\ell')]
=\delta_{\ell+\ell',0}\frac{
(\gamma\gamma')^{\frac{|\ell|}{2}}(rs)^{-\frac{\ell (\alpha_i|\alpha_j)}2}[\,\ell\,
(\alpha_i|\alpha_j)\,]}{|\ell|}
\cdot\frac{\gamma^{|\ell|}-\gamma'^{|\ell|}}{r-s}.\\
({\rm D3})&&[a_i(\ell),\om_j^{{\pm }1}]=[a_i(\ell),\om_j'^{\pm
1}]=0.\\
%\end{eqnarray*}
%\begin{eqnarray*}
({\rm D4})&& D\,x_i^{\pm}(k)\,D^{-1}=r_i^k\, x_i^{\pm}(k), \qquad
D'\,x_i^{\pm}(k)\,D'^{\,-1}=s_i^k\, x_i^{\pm}(k), \\
&&D\,
a_i(\ell)\,D^{-1}=r_i^l\,a_i(\ell), \qquad\qquad D'\,
a_i(\ell)\,D'^{\,-1}=s_i^l\,a_i(\ell).\\
({\rm D5}) &&\om_i\,x_j^{\pm}(k)\, \om_i^{-1} =  \langle j,
i\rangle^{\pm 1} x_j^{\pm}(k), \qquad \om'_i\,x_j^{\pm}(k)\,
\om_i'^{\,-1} =  \langle i, j\rangle ^{\mp1}x_j^{\pm}(k).\\
({\rm D6_1})&&[\,a_i(\ell),x_j^{\pm}(k)\,]=\pm\frac{
(rs)^{\frac{\ell(c-(\alpha_i|\alpha_j))}{2}}[\,\ell
(\alpha_i|\alpha_j)\,]}{\ell}\gamma'^{\pm\frac{\ell}2}x_j^{\pm}(\ell{+}k),
 \quad \textit{for } \ \ell>0, \\
({\rm D6_2})&&[\,a_i(\ell),x_j^{\pm}(k)\,]=\pm\frac{
(rs)^{\frac{-\ell (c+(\alpha_i|\alpha_j))}{2}}[\,\ell\,
(\alpha_i|\alpha_j)\,]}{\ell}\gamma^{\pm\frac{\ell}2}x_j^{\pm}(\ell{+}k),
\quad \textit{for } \ell<0. \\
({\rm D7})&&x_i^{\pm}(k{+}1)\,x_j^{\pm}(k') - \la j,i\ra^{\pm1} x_j^{\pm}(k')\,x_i^{\pm}(k{+}1)\\
&&=-\Bigl(\la j,i\ra\la
i,j\ra^{-1}\Bigr)^{\pm\frac1{2}}\,\Bigl(x_j^{\pm}(k'{+}1)\,x_i^{\pm}(k)-\la
i,j\ra^{\pm1} x_i^{\pm}(k)\,x_j^{\pm}(k'{+}1)\Bigr),\\
({\rm
D8})&&[\,x_i^{+}(k),~x_j^-(k')\,]=\frac{\delta_{ij}}{r_i-s_i}\Big(\gamma'^{-k}\,{\gamma}^{-\frac{k+k'}{2}}\,
\psi_i(k{+}k')-\gamma^{k'}\,\gamma'^{\frac{k+k'}{2}}\,\varphi_i(k{+}k')\Big),
\end{eqnarray*}
where
$\psi_i(m)$, $\varphi_i(-m)~(m\in \mathbb{Z}_{\geq 0})$ with
$\psi_i(0)=\om_i$ and $\varphi_i(0)=\om_i'$ are defined by:

$$\sum\limits_{m=0}^{\infty}\psi_i(m) z^{-m}=\om_i \exp \Big(
(r{-}s)\sum\limits_{\ell=1}^{\infty}
 a_i(\ell)z^{-\ell}\Big),\quad \bigl(\psi_i(-m)=0, \ \forall\;m>0\bigr);
$$
$$\sum\limits_{m=0}^{\infty}\varphi_i(-m) z^{m}=\om'_i \exp
\Big({-}(r{-}s)
\sum\limits_{\ell=1}^{\infty}a_i(-\ell)z^{\ell}\Big), \quad
\bigl(\varphi_i(m)=0, \ \forall\;m>0\bigr).
$$

\medskip
\noindent$({\rm D9_1})$ \ $~~~Sym_{n_1, n_2}\Big(x_i^{\pm}(n_1)
x_i^{\pm}(n_2)x_j^{\pm}(k)-(r_i^{\pm1}{+}s_i^{\pm1})\,x_i^{\pm}(n_1)x_j^{\pm}(k)x_i^{\pm}(n_2)$

$\hskip3cm +(r_is_i)^{\pm1}
x_j^{\pm}(k)x_i^{\pm}(n_1)x_i^{\pm}(n_2)\Big)=0, $

$\hskip 6.5cm  \textit{for }
\ a_{ij}=-1 \ \textit{and}\  1\leqslant i < j \leqslant n;$

\medskip
\noindent $({\rm D9_2})$ \ $~~~Sym_{n_1, n_2}\Big(x_i^{\pm}(n_1)
x_i^{\pm}(n_2)x_j^{\pm}(k)-(r_i^{\mp1}{+}s_i^{\mp1})\,x_i^{\pm}(n_1)x_j^{\pm}(k)x_i^{\pm}(n_2)$

$\hskip3cm +(r_is_i)^{\mp1}
x_j^{\pm}(k)x_i^{\pm}(n_1)x_i^{\pm}(n_2)\Big)=0, $

$\hskip 6.5cm  \textit{for }
\ a_{ij}=-1\  \textit{and}\  1\leqslant j < i \leqslant n;$

\medskip
\noindent $({\rm D9_3})$ \ $Sym_{n_1,
n_2,n_3}\Big(x_i^{\pm}(n_1)
x_i^{\pm}(n_2)x_i^{\pm}(n_3)x_j^{\pm}(k)
$\\

$-\,
(r_i^{\pm2}{+}(r_is_i)^{\pm1}(+)s_i^{\pm2})\,x_i^{\pm}(n_1)x_i^{\pm}(n_2)x_j^{\pm}(k)x_i^{\pm}(n_3)
$\\

$+\,(r_is_i)^{\pm1}(r_i^{\pm2}{+}(r_is_i)^{\pm1}(+)s_i^{\pm2})
x_i^{\pm}(n_1)x_j^{\pm}(k)x_i^{\pm}(n_2)x_i^{\pm}(n_3)
$\\

$-\,(r_is_i)^{\pm3}x_j^{\pm}(k)x_i^{\pm}(n_1)x_i^{\pm}(n_2)x_i^{\pm}(n_3)\Big)=0,
\quad \textit{for } \ a_{ij}=-2$.
\end{defi}
%\begin{remark}\, Note that the center $c$ in the relations
%$({\rm D2}),({\rm D6})$ of the definition, which is relevant with the level of the representation of ${\mathcal U}_{r,s}(C_n^{(1)})$ (for detailed reasons, see \cite{HZ}).
%When $r=s^{-1}=q$, this is the Drinfeld realization of one-parameter quantum affine algebras.
%\end{remark}
\begin{remark}  Notice that the defining relations
$(\textrm{D7})$, $(\textrm{D8})$ can be written equivalently by
virtue of generating functions (see \cite{HZ}) as follows:
\begin{gather*}
(z-(\lg i,j\rg\lg j,i\rg)^{\pm\frac1{2}}w)\,x_i^{\pm}(z)x_j^{\pm}(w)
=(\lg j,i\rg^{\pm1}z-(\lg j,i\rg\lg
i,j\rg^{-1})^{\pm\frac1{2}}w)\,x_j^{\pm}(w)\,x_i^{\pm}(z),\tag{\textrm{D7}$'$}\\
[\,x_i^+(z),
x_j^-(w)\,]=\frac{\delta_{ij}}{(r_i-s_i)zw}\Big(\delta(\frac{w{\gamma'}^{-1}}{z})
\psi_i(w\gamma^{\frac{1}2})
-\delta(\frac{w\gamma^{-1}}{z})\varphi_i(w\gamma'^{-\frac1{2}})\Big),\tag{\textrm{D8}$'$}
\end{gather*}
where $\delta(z)=\sum_{n\in\mathbb{Z}}z^n$,
$x_i^{\pm}(z)=\sum\limits_{n\in\mathbb{Z}}x_i^{\pm}(n)z^{-n-1}$,
$\psi_i(z) = \sum_{m \in \mathbb{Z}_+}\psi_i(m) z^{-m}$, and
$\varphi_i(z) = \sum_{n \in -\mathbb{Z}_+}\varphi_i(n) z^{-n}$.
\end{remark}

\medskip
As one of crucial observations of considering the compatibilities of
the defining system above, we have
\begin{prop}
There exists the $\mathbb{Q}$-algebra antiautomorphism $\tau$ of
${\mathcal U}_{r,s}(C_n^{(1)})$  such that
$\tau(r)=s$, $\tau(s)=r$,
$\tau(\la\om_i',\om_j\ra^{\pm1})=\la\om_j',\om_i\ra^{\mp1}$ and
\begin{gather*}
\tau(\om_i)=\om_i',\quad \tau(\om_i')=\om_i,\quad
\tau(\gamma)=\gamma',\quad \tau(\gamma')=\gamma,\quad\tau(D)=D',\quad \tau(D')=D,\\
\tau(x_i^{\pm}(m))=x_i^{\mp}(-m), \quad
\tau(a_i(\ell))=a_i(-\ell),\\
\tau(\psi_i(m))=\varphi_i(-m), \quad\tau(\varphi_i(-m))=\psi_i(m),
\end{gather*}
and $\tau$ preserves each defining relation $($\hbox{{\rm D}n}$)$ in
Definition 3.1 for $n=1, \cdots, 9$.\hfill\qed
\end{prop}
\begin{remark}The defining relations $(\rm{D1})$---$(\rm{D9})$ ensure that
$\mathcal{U}_{r,s}(\mathrm{C}_n^{(1)})$ has a triangular
decomposition:
$$\mathcal{U}_{r,s}(\mathrm{C}_n^{(1)})=
\mathcal{U}_{r,s}(\widetilde{\frak{n}}^-)\bigotimes\mathcal{U}_{r,s}^0(\mathrm{C}_n^{(1)})
\bigotimes\mathcal{U}_{r,s}(\widetilde{\frak{n}}),$$ where
$\mathcal{U}_{r,s}(\widetilde{\frak{n}}^\pm)=\bigoplus_{\alpha\in\dot
Q^\pm}\mathcal{U}_{r,s}(\widetilde{\frak{n}}^\pm)_\alpha$ is
generated respectively by $x_i^\pm(k)$ $(i\in I)$, and
$\mathcal{U}_{r,s}^0(\mathrm{C}_n^{(1)})$ is the subalgebra
generated by $\om_i^{\pm1}$, $\om_i'^{\pm1}$,
$\gamma^{\pm\frac1{2}}$, $\gamma'^{\pm\frac1{2}}$, $D^{\pm1}$,
$D'^{\pm1}$ and $a_i(\pm\ell)$ for $i\in I$, $\ell\in \mathbb{N}$.
Namely, $\mathcal{U}_{r,s}^0(\mathrm{C}_n^{(1)})$ is generated by
the toral subalgebra $\mathcal{U}_{r,s}(\mathrm{C}_n^{(1)})^0$ and
the quantum Heisenberg subalgebra $\mathcal
H_{r,s}(\mathrm{C}_n^{(1)})$ generated by those quantum imaginary
root vectors $a_i(\pm\ell)$ $(i\in I$, $\ell\in \mathbb{N})$.
\end{remark}

\subsection{Drinfeld Isomorphism}  \ To obtain the
isomorphism between the above two realizations for the two-parameter
quantum affine algebra $U_{r,s}(C_n^{(1)})$, we need the following
notations and definitions.
\begin{defi}\ $($\cite{HRZ}$)$ The quantum Lie brackets $[\,a_1,\cdots,
a_s\,]_{(q_1,\,\cdots,\, q_{s-1})}$ and $[\,a_1, a_2, \cdots,
a_s\,]_{\la q_1,\,q_2,\,\cdots, \,q_{s-1}\ra}$ are defined
recursively by
\begin{gather*}    [\,a_1, a_2\,]_{v_1}=a_1a_2-v_1\,a_2a_1,\\
   [\,a_1, a_2, \cdots, a_s\,]_{(v_1,\,v_2,\,\cdots,
  \,v_{s-1})}=[\,a_1, \cdots, [a_{s-1},
  a_s\,]_{v_1}]_{(v_2,\,\cdots,\,v_{s-1})},\\
  [\,a_1, a_2, \cdots, a_s\,]_{\la v_1,\,v_2,\,\cdots,
  \,v_{s-1}\ra}=[\,[\,a_1, a_2]_{v_1} \cdots, a_{s-1}\,]_{\la
  v_2,\,\cdots,\,v_{s-2}\ra},
 \end{gather*}

\noindent
for $q_i\in \mathbb K^*=\mathbb{K}\backslash \{0\}$.
\end{defi}

For $q\in\mathbb K^*$, the following identities follow from the
definition
\begin{eqnarray}
&&[\,a, bc\,]_v=[\,a, b\,]_q\,c+q\,b\,[\,a, c\,]_{\frac{v}q},\\
&&[\,ab, c\,]_v=a\,[\,b, c\,]_q+q\,[\,a, c\,]_{\frac{v}q}\,b, \\
&&{\label{b:3}[\,a,[\,b,c\,]_u\,]_v=[\,[\,a,b\,]_q,
c\,]_{\frac{uv}q}+q\,[\,b,[\,a,c\,] _{\frac{v}q}\,]_{\frac{u}q},}\\
&&{\label{b:4}[\,[\,a,b\,]_u,c\,]_v=[\,a,[\,b,c\,]_q\,]_{\frac{uv}q}+q\,[\,[\,a,c\,]
_{\frac{v}q},b\,]_{\frac{u}q}}.
\end{eqnarray}

In particular, we have (see \cite{HRZ}, \cite{HZ})
\begin{eqnarray}
&&{\label{b:5}[\,a, [\,b_1, \cdots, b_s\,]_{(v_1,\,\cdots,\,
v_{s-1})}\,]=\sum\limits_i[\,b_1,\cdots,[\,a, b_i\,],
\cdots,b_s\,]_{(v_1,\,\cdots,\, v_{s-1})},\hskip0.3cm }\\
&&{\label{b:6}[\,a, a, b\,]_{(u,\,
v)}=a^2b-(u{+}v)\,aba+(uv)\,ba^2=(uv)[\,b, a, a\,]_{\la u^{-1},v^{-1}\ra},} \\
&&[\,a, a, a, b\,]_{(u^2,\,uv,\,v^2)}=a^3b-[3]_{u,v}\,a^2ba
+(uv)[3]_{u,v}aba^2-(uv)^3ba^3.\nonumber
\end{eqnarray}
where $[n]_{u,v}=\frac{u^n{-}v^n}{u{-}v}$,
$[n]_{u,v}!:=[n]_{u,v}\cdots [2]_{u,v}[1]_{u,v}$, $\left[n\atop
m\right]_{u,v}:=\frac{[n]_{u,v}!}{[m]_{u,v}![n-m]_{u,v}!}$.

\begin{defi} $($\cite{HRZ}$)$
For $\alpha,\,\beta\in\dot Q^+$ $($the positive root lattice of $C_n)$,
and $x_\alpha^\pm(k),\,x_\beta^\pm(k')$ $\in
\mathcal{U}_{r,s}(\widetilde{\frak n}^\pm)$, define the {\rm affine}
quantum Lie bracket as follows:
$$\bigl[\,x_\alpha^\pm(k),\,
x_\beta^\pm(k')\,\bigr]_{\la\om'_\alpha,\om_\beta\ra^{\mp1}}:=
x_\alpha^\pm(k)\,x_\beta^\pm(k')-\la\om'_\alpha,\om_\beta\ra^{\mp1}
x_\beta^\pm(k')\,x_\alpha^\pm(k).
$$
\end{defi}
By the definition above, formula $(\textrm{D7})$ will take the
convenient form as
\begin{eqnarray*}
[x_i^\pm(k), x_j^\pm(k'{+}1)]_{\la i, j\ra^{\mp1}}=-(\la j, i\ra\la
i, j\ra^{-1})^{\pm\frac1{2}}[x_j^\pm(k'), x_i^\pm(k{+}1)]_{\la j,
i\ra^{\mp1}}.
\end{eqnarray*}
Especially, we have
\begin{gather}
{\label{b:7}[\,x^{-}_2(0),\, x_1^{-}(1)\,]_{s^{\frac{1}{2}}}
=-(rs)^{\frac{1}{4}}[\,x^{-}_1(0),\, x_2^{-}(1)\,]_{r^{-\frac{1}{2}}},}\\
{\label{b:8}[\,x_i^\pm(k), x_i^\pm(k{+}1)\,]_{\la i, i\ra^{\mp1}}=0.}
\end{gather}
By (3.6) \& (3.8), the $(r,s)$-Serre relations (D$9_1$), (D$9_2$) \& (D$9_3$)
for $n_i=n_j=\ell$ in the case of $a_{ij}=-1$ and $a_{n-1\, n}=-2$ can be
reformulated respectively as:
\begin{eqnarray}
&&{\label{b:9}[\,x_i^{\pm}(\ell),
\,x_i^{\pm}(\ell),\,
x_{i-1}^{\pm}(k)\,]_{(s_i^{\pm1},\,r_i^{\pm1})}=0, } \\
&&{\label{b:10}[\,x_i^{\pm}(\ell),
\,x_i^{\pm}(\ell),\,
x_{i+1}^{\pm}(k)\,]_{(s_i^{\mp1},\,r_i^{\mp1})}=0,}  \\
&&{\label{b:11}[\,x_{n-1}^{\pm}(\ell),\, x_{n-1}^{\pm}(\ell),\,
x_{n-1}^{\pm}(\ell),\,x_n^{\pm}(k)\,]_{(r_{n-1}^{\pm2},\,s_{n-1}^{\pm1}s_{n-1}^{\pm1},
\,s_{n-1}^{\pm2})}=0.}
\end{eqnarray}

For $2\leqslant i \leqslant n-1$, let us set some notations for later use,
\begin{equation}
{\label{b:12}
\begin{split}
x_{\alpha_{1,i}}^-(1)=[\,x_i^-(0),\,\cdots,\,x_2^-(0),\,x_1^-(1)
\,]_{(s^{\frac{1}{2}},\,\cdots,\,s^{\frac{1}{2}})}
\end{split}}
\end{equation}
\medskip

\begin{equation}
{\label{b:13}
\begin{split}
x_{\alpha_{1,n}}^-(1)=[\,x_n^-(0),\,\cdots,\,x_2^-(0),\,x_1^-(1)
\,]_{(s^{\frac{1}{2}},\,\cdots,\,s^{\frac{1}{2}}, s)}
\end{split}}
\end{equation}
\medskip

\begin{equation}
{\label{b:14}
\begin{split}
&x_{\beta_{1,i}}^-(1)\\
&\quad=[x_{i}^-(0), x_{i+1}^-(0), \cdots, x_{n}^-(0), \cdots, x_2^-(0), x_1^-(1)]_{(s^{\frac{1}{2}}, \cdots, s^{\frac{1}{2}}, s, r^{-\frac{1}{2}}, \cdots, r^{-\frac{1}{2}})},
\end{split}}
\end{equation}

\begin{equation}
{\label{b:15}
\begin{split}
&x_{\beta_{1,1}}^-(1)\\
&\quad=[x_{1}^-(0), x_{2}^-(0), \cdots, x_{n}^-(0), \cdots, x_2^-(0), x_1^-(1)]_{(s^{\frac{1}{2}}, \cdots, s^{\frac{1}{2}}, s, r^{-\frac{1}{2}}, \cdots, r^{-\frac{1}{2}}, (rs)^{-\frac{1}{2}})}.
\end{split}}
\end{equation}

\begin{remark} In particular, we denote that,
\begin{equation*}
\begin{split}
&x_{\theta}^+(-1):=x_{\beta_{1,1}}^+(-1)\\
=&\Big[\,x_1^+(-1), x_2^+(0), \cdots, x_{n}^+(0), \cdots, x_1^+(0)\Big]_{\langle r^{\frac{1}{2}}, \cdots, r^{\frac{1}{2}}, r,
s^{-\frac{1}{2}},\cdots,\,s^{-\frac{1}{2}},(rs)^{-\frac{1}{2}}\rangle} ,\\
&x_{\theta}^-(1):=x_{\beta_{1\,1}}^-(1)\\
=&\Big[\,x_1^-(0),\,\cdots,\,x_n^-(0),\,\cdots,\,x_2^-(0),\,x_1^-(1)
\,\Big]_{(s^{\frac{1}{2}},\,\cdots,\,s^{\frac{1}{2}},s,\,
r^{-\frac{1}{2}},\,\cdots,\,r^{-\frac{1}{2}},(rs)^{-\frac{1}{2}})}.
\end{split}
\end{equation*}
\end{remark}

The following two lemmas will be used in the proof of the main
theorem of this section.
\begin{lemm} Using the above notations, one has,
\begin{eqnarray}
&&{\label{b:16} [\,x_i^-{(0)},\, x_{\alpha_{1\,i}}^-(1)\,]_{r^{\frac{1}{2}}}=0, \, \hbox{for}\  2\leqslant i \leqslant n-1 ;}\\
&&{\label{b:17} [\,x_{i}^-{(0)},\, x_{\alpha_{1\,i+1}}^-(1)\,]=0, \, \hbox{for}\  2\leqslant i \leqslant n.}
\end{eqnarray}
\end{lemm}

\begin{proof} (\ref{b:16}): Using the above notations and (\ref{b:3}), one has
\begin{eqnarray*}
&&[\,x_i^-{(0)},\, x_{\alpha_{1\,i}}^-(1)\,]_{r^{\frac{1}{2}}}
\quad {\hbox{(by definition)}}\\
&=&[\,x_i^-(0),\,[\,x_{i}^-(0),\,x_{i-1}^-(0),\,x_{\alpha_{1,i-2}}^-(1)
\,]_{(s^{\frac{1}{2}},\,s^{\frac{1}{2}})}\,]_{r^{\frac{1}{2}}}
\quad {\hbox{(using (\ref{b:3}))}}\\
&=&[\,x_i^-(0),\,[\,x_{i}^-(0),\,x_{i-1}^-(0)\,]_{s^{\frac{1}{2}}}
,\,x_{\alpha_{1\,i-2}}^-(1)\,]_{(s^{\frac{1}{2}},\,r^{\frac{1}{2}})}
\quad {\hbox{(using (\ref{b:3}))}}\\
&&+s^{\frac{1}{2}}[\,x_i^-(0),\,x_{i-1}^-(0),\,
\underbrace{[\,x_{i}^-(0),\,x_{\alpha_{1\,i-2}}^-(1)\,]}
\,]_{(1,\,r^{\frac{1}{2}})}\quad {\hbox{(=0 by (\ref{b:5}) )}}\\
&=&[\,\underbrace{[\,x_i^-(0),\,[\,x_{i}^-(0),\,x_{i-1}^-(0)\,]_{s^{\frac{1}{2}}}
\,]_{r^{\frac{1}{2}}}}
,\,x_{\alpha_{1\,i-2}}^-(1)\,]_{s^{\frac{1}{2}}}\quad {\hbox{(=0 by $(\ref{b:9})$)}}\\
&&+r^{\frac{1}{2}}[\,[\,x_{i}^-(0),\,x_{i-1}^-(0)\,]_{s^{\frac{1}{2}}}
,\,\underbrace{[\,x_{i}^-(0),\,x_{\alpha_{1\,i-2}}^-(1)\,]}
\,]_{r^{-\frac{1}{2}}s^{\frac{1}{2}}}\quad {\hbox{(=0 by (\ref{b:5}))}}\\
&=&0.
\end{eqnarray*}
(3.17): First note that
\begin{eqnarray*}
\begin{split}
&[\,x_i^-{(0)},\, x_{\alpha_{1\,i+1}}^-(1)\,]_{r^{-\frac{1}{2}}s^{\frac{1}{2}}}
\quad {\hbox{(by definition)}}\\
&\quad=[\,x_i^-(0),\,[\,x_{i+1}^-(0),\,x_i^-(0),\,x_{\alpha_{1\,i-1}}^-(1)
\,]_{(s^{\frac{1}{2}},\,s^{\frac{1}{2}})}\,]_{r^{-\frac{1}{2}}s^{\frac{1}{2}}}
\quad {\hbox{(using (\ref{b:3}))}}\\
&\quad=[\,x_i^-(0),\,[\,x_{i+1}^-(0),\,x_i^-(0)\,]_{s^{\frac{1}{2}}}
,\,x_{\alpha_{1\,i-1}}^-(1)\,]_{(s^{\frac{1}{2}},\,r^{-\frac{1}{2}}s^{\frac{1}{2}})}
\quad {\hbox{(using (\ref{b:3}))}}\\
&\quad+s^{\frac{1}{2}}[\,x_i^-(0),\,x_i^-(0),\,
\underbrace{[\,x_{i+1}^-(0),\,x_{\alpha_{1\,i-1}}^-(1)\,]}
\,]_{(1,\,r^{-\frac{1}{2}}s^{\frac{1}{2}})}\quad {\hbox{(=0 by (\ref{b:5}))}}\\
\end{split}
\end{eqnarray*}
\begin{eqnarray*}
\begin{split}
&\quad=[\,\underbrace{[\,x_i^-(0),\,[\,x_{i+1}^-(0),\,x_i^-(0)\,]_{s^{\frac{1}{2}}}
\,]_{r^{-\frac{1}{2}}}}
,\,x_{\alpha_{1\,i-1}}^-(1)\,]_{s}\quad {\hbox{(=0 by (\ref{b:10}))}}\\
&\quad+r^{-\frac{1}{2}}[\,[\,x_{i+1}^-(0),\,x_i^-(0)\,]_{s^{\frac{1}{2}}}
,\,[\,x_{i}^-(0),\,x_{\alpha_{1\,i-1}}^-(1)\,]_{s^{\frac{1}{2}}}
\,]_{(rs)^{\frac{1}{2}}}\quad {\hbox{(by definition \& (\ref{b:4}))}}\\
&\quad=r^{-\frac{1}{2}}[\,x_{i+1}^-(0),\,
\underbrace{[\,x_i^-(0),\,x_{\alpha_{1\,i}}^-(1)\,]_{r^{\frac{1}{2}}}}
\,]_{s}\quad {\hbox{(=0 by (\ref{b:16}))}}\\
&\quad+[\,[\,x_{i+1}^-(0),\,x_{\alpha_{1\,i}}^-(1)\,]_{s^{\frac{1}{2}}}
,\,x_{i}^-(0)\,]_{r^{-\frac{1}{2}}s^{\frac{1}{2}}}
\quad {\hbox{(by definition)}}\\
&\quad=[\,x_{\alpha_{1\,i+1}}^-(1),\,x_{i}^-(0)\,]_{r^{-\frac{1}{2}}s^{\frac{1}{2}}},
\end{split}
\end{eqnarray*}
which implies that
$(1+r^{-\frac{1}{2}}s^{\frac{1}{2}})[\,x_i^-(0),\,
x_{\alpha_{1\,i+1}}^-(1)\,]=0$, then we get the required relation.
\end{proof}

\begin{lemm}{\label{d:16} The following relations hold:
\begin{eqnarray}
&&{\label{b:18}[\,x_1^-{(0)},\, x_{\beta_{1\,3}}^-(1)\,]_{s^{-\frac{1}{2}}}=0; }\\
&&{\label{b:19}[\,x_i^-{(0)},\, x_{\beta_{1\,i+2}}^-(1)\,]=0; \qquad\quad \hbox{for} \qquad 2\leqslant i \leqslant n-2;}\\
&&{\label{b:20}[\,x_i^-{(0)},\, x_{\beta_{1\,i}}^-(1)\,]_{s^{-\frac{1}{2}}}=0, \qquad\quad \hbox{for} \qquad  2\leqslant i \leqslant n-1;} \\
&&{\label{b:21}[\,x_n^-{(0)},\, x_{\alpha_{1\,n}}^-(1)\,]_r=0.}
\end{eqnarray}}
\end{lemm}
\begin{proof} (\ref{b:18}): It is easy to obtain that
\begin{eqnarray*}
\begin{split}
&[\,x_1^-{(0)},\, x_{\beta_{1\,3}}^-(1)\,]_{s^{-\frac{1}{2}}}
\quad {\hbox{(by definition)}}\\
&\quad=[x_1^-(0), [x_{3}^-(0), \cdots,
x_{n}^-(0), \cdots, x_3^-(0), x_2^-(0), x_{1}^-(1)]\\
&\quad\hskip1.5cm _{(s^{\frac{1}{2}}, \cdots, s^{\frac{1}{2}}, s
, r^{-\frac{1}{2}}, \cdots, r^{-\frac{1}{2}}, (rs)^{-\frac{1}{2}})}\,]_{s^{-\frac{1}{2}}}\quad {\hbox{(using (\ref{b:3}) )}}\\
&\quad=[\,x_{3}^-(0),\,\cdots, x_{n}^-(0),\,\cdots,x_3^-(0),\,
[\,x_1^-(0),\,\underbrace{[\,x_2^-(0),\,x_{1}^-(1)\,]_{s^{\frac{1}{2}}}}
\,]_{s^{-\frac{1}{2}}}\\
&\quad\hskip2.5cm\,]_{(s^{\frac{1}{2}},\cdots,s^{\frac{1}{2}},\,s
,\,r^{-\frac{1}{2}},\cdots,r^{-\frac{1}{2}},\,(rs)^{-\frac{1}{2}})}
\,]_{s^{-\frac{1}{2}}}\qquad {\hbox{(using (\ref{b:3}))}}\\
&\quad=-(rs)^{\frac{1}{4}}[\,x_{3}^-(0),\,\cdots, x_{n}^-(0),\,\cdots,x_3^-(0),\,
\underbrace{[\,x_1^-(0),\,[\,x_1^-(0),\,x_{2}^-(1)\,]_{r^{-\frac{1}{2}}}
\,]_{s^{-\frac{1}{2}}}}\\
&\quad\hskip2.5cm\,]_{(s^{\frac{1}{2}},\cdots,s^{\frac{1}{2}},\,s
,\,r^{-\frac{1}{2}},\cdots,r^{-\frac{1}{2}},\,(rs)^{-\frac{1}{2}})}
\,]_{s^{-\frac{1}{2}}}\qquad {\hbox{(=0 by (\ref{b:10}))}}\\
&\quad=0.
\end{split}
\end{eqnarray*}
(\ref{b:19}): One has
\begin{eqnarray*}
\begin{split}
&[\,x_i^-{(0)},\, x_{\beta_{1\,i+2}}^-(1)\,]
\qquad {\hbox{(by definition)}}\\
&\quad=[\,x_i^-(0),\,[\,x_{i+2}^-(0),\,\cdots,
x_{n}^-(0),\,\cdots,x_{i+2}^-(0),\,
x_{\alpha_{1\,i+1}}^-(1)\,]\\
&\quad\hskip2.4cm \cdots]_{(s^{\frac{1}{2}},\cdots,s^{\frac{1}{2}},\,s
,\,r^{-\frac{1}{2}},\cdots,r^{-\frac{1}{2}},\,(rs)^{-\frac{1}{2}})}
\,]\quad {\hbox{(using (\ref{b:3}))}}\\
&\quad=[\,x_{i+2}^-(0),\,\cdots, x_{n}^-(0),\,\cdots,x_{i+2}^-(0),\,
\underbrace{[\,x_i^-(0),\,x_{\alpha_{1\,i+1}}^-(1)\,]}
\,]\\
&\quad\hskip2.4cm \cdots]_{(s^{\frac{1}{2}},\cdots,s^{\frac{1}{2}},\,s
,\,r^{-\frac{1}{2}},\cdots,r^{-\frac{1}{2}},\,(rs)^{-\frac{1}{2}})}\qquad {\hbox{(=0 by (\ref{b:17}))}}\\
&\quad=0.
\end{split}
\end{eqnarray*}
(\ref{b:20}): We can get without difficulty:
\begin{eqnarray*}
\begin{split}
&[\,x_i^-{(0)},\, x_{\beta_{1\,i}}^-(1)\,]_{s^{-\frac{1}{2}}}
\quad {\hbox{(by definition)}}\\
&\quad=[\,x_i^-(0),\,[\,x_{i}^-(0),\,x_{i+1}^-(0),\,x_{\beta_{1\,i+2}}^-(1)
\,]_{(r^{-\frac{1}{2}},\,r^{-\frac{1}{2}})}\,]_{s^{-\frac{1}{2}}}
\quad {\hbox{(using (\ref{b:3}))}}\\
&\quad=[\,x_i^-(0),\,[\,x_{i}^-(0),\,x_{i+1}^-(0)\,]_{r^{-\frac{1}{2}}}
,\,x_{\beta_{1\,i+2}}^-(1)\,]_{(r^{-\frac{1}{2}},\,s^{-\frac{1}{2}})}
\quad {\hbox{(using (\ref{b:3}))}}\\
&\quad+r^{-\frac{1}{2}}[\,x_i^-(0),\,x_{i+1}^-(0),\,
\underbrace{[\,x_{i}^-(0),\,x_{\beta_{1\,i+2}}^-(1)\,]}
\,]_{(1,\,s^{-\frac{1}{2}})}\quad {\hbox{(=0 by (\ref{b:19}))}}\\
&\quad=[\,\underbrace{[\,x_i^-(0),\,[\,x_{i}^-(0),\,x_{i+1}^-(0)\,]_{r^{-\frac{1}{2}}}
\,]_{s^{-\frac{1}{2}}}}
,\,x_{\beta_{1\,i+2}}^-(1)\,]_{r^{-\frac{1}{2}}}\quad {\hbox{(=0 by (\ref{b:10}))}}\\
&\quad+s^{-\frac{1}{2}}[\,[\,x_{i}^-(0),\,x_{i-1}^-(0)\,]_{r^{-\frac{1}{2}}}
,\,\underbrace{[\,x_{i}^-(0),\,x_{\beta_{1\,i+2}}^-(1)\,]}
\,]_{r^{-\frac{1}{2}}s^{\frac{1}{2}}}\quad {\hbox{(=0 by (\ref{b:19}))}}\\
&\quad=0.
\end{split}
\end{eqnarray*}
(\ref{b:21}):  By calculating, we obtain
\begin{eqnarray*}
\begin{split}
&[\,x_n^-{(0)},\, x_{\alpha_{1\,n}}^-(1)\,]_{r}
\quad {\hbox{(by definition)}}\\
&\quad=[\,x_n^-(0),\,[\,x_{n}^-(0),\,x_{n-1}^-(0),\,x_{\alpha_{1\,n-2}}^-(1)
\,]_{(s^{\frac{1}{2}},\,s)}\,]_{r}
\quad {\hbox{(using (\ref{b:3}))}}\\
&\quad=[\,x_n^-(0),\,[\,x_{n}^-(0),\,x_{n-1}^-(0)\,]_{s}
,\,x_{\alpha_{1\,n-2}}^-(1)\,]_{(s^{\frac{1}{2}},\,r)}
\quad {\hbox{(using (\ref{b:3}))}}\\
&\quad+s[\,x_n^-(0),\,x_{n-1}^-(0),\,
\underbrace{[\,x_{n}^-(0),\,x_{\alpha_{1\,n-2}}^-(1)\,]}
\,]_{(s^{-\frac{1}{2}},\,r)}\quad {\hbox{(=0 by (\ref{b:5}) and $(\textrm{D9}_1)$)}}\\
&\quad=[\,\underbrace{[\,x_n^-(0),\,[\,x_{n}^-(0),\,x_{n-1}^-(0)\,]_{s}
\,]_{r}}
,\,x_{\alpha_{1,n-2}}^-(1)\,]_{s^{\frac{1}{2}}}\quad {\hbox{(=0 by (\ref{b:9}))}}\\
&\quad+r[\,[\,x_{n}^-(0),\,x_{n-1}^-(0)\,]_{s}
,\,\underbrace{[\,x_{n}^-(0),\,x_{\alpha_{1\,n-2}}^-(1)\,]}
\,]_{r^{-\frac{1}{2}}s^{\frac{1}{2}}}\quad {\hbox{(=0 by (\ref{b:5}))}}\\
&\quad=0.
\end{split}
\end{eqnarray*}
\end{proof}

Now we turn to give one of our main theorems as follows.

\begin{theo} \ For non-simply-laced Lie algebra $C_n^{(1)}$,
 there exists an algebra isomorphism
$\Psi: U_{r,s}(C_n^{(1)}) \longrightarrow {\mathcal
U}_{r,s}(C_n^{(1)})$ defined by: for $i\in I,$
\begin{eqnarray*}
\omega_i&\longmapsto& \om_i\\
\omega'_i&\longmapsto& \om'_i \\
\omega_0&\longmapsto& \gamma'^{-1}\, \om_{\theta}^{-1}\\
\omega'_0&\longmapsto& \gamma^{-1}\, \om_{\theta}'^{-1}\\
\gamma^{\pm\frac{1}2}&\longmapsto& \gamma^{\pm\frac{1}2}\\
\gamma'^{\,\pm\frac{1}2}&\longmapsto& \gamma'^{\,\pm\frac{1}2}\\
D^{\pm1}&\longmapsto& D^{\pm1}\\
D'^{\,\pm1}&\longmapsto& D'^{\,\pm1}\\
e_i&\longmapsto& x_i^+(0)\\
f_i&\longmapsto& x_i^-(0)\\
e_0&\longmapsto& a\,
x^-_{\theta}(1)\cdot(\gamma'^{-1}\,\om_{\theta}^{-1})\\
f_0 &\longmapsto& \tau\Bigl(a\,x^-_{\theta}(1)
      \,\cdot(\gamma'^{-1}\,\om_{\theta}^{-1})\Bigr)
      =a\,(\gamma^{-1}\,{\om'}_{\theta}^{-1})\cdot
      x_{\theta}^+(-1)
\end{eqnarray*}
where $\om_{\theta}=\om_{i_1}\,\cdots\,
\om_{i_{h-1}},\,\om'_{\theta}=\om'_{i_1}\,\cdots\, \om'_{i_{h-1}}$,
and $a=(rs)^{\frac{n-2}{2}}[2]_1^{-1}$.
\end{theo}
\begin{remark}  We note that $\tau(a)=a$.
\end{remark}

\subsection{The proof of Theorem 3.10} \  In this subsection, we
prove Theorem 3.10. Let $E_i,\,F_i$ ($i\in I_0$) and $\om_0$,
$\om_0'$ denote the images of $e_i,\,f_i$ ($i\in I_0$) and $\om_0$,
$\om_0'$ in the algebra ${\mathcal U}_{r,s}(C_n^{(1)})$,
respectively.

Denote by $\mathcal{U}'_{r,s}(C_n^{(1)})$ the subalgebra of
$\mathcal{U}_{r,s}(C_n^{(1)})$ generated by
$E_i,\,F_i,\,\om_i^{\pm1}$, $\om_i'^{\pm1}$ ($i\in I_0$),
$\gamma^{\pm\frac{1}2},\, \gamma'^{\,\pm\frac{1}2}$, $D^{\pm1}$ and
$D'^{\pm1}$, that is,
$$
{\mathcal U}'_{r,s}(C_n^{(1)}):=\left.\left\langle\, E_i,\, F_i,\,
\om_i^{\pm1},\, {\om}_i'^{\,\pm1}\;,\, \gamma^{\pm\frac{1}2},\,
\gamma'^{\,\pm\frac{1}2},\, D^{\pm1},\, D'^{\,\pm1}\; \right| \;i\in
I_0\;\right\rangle.
$$

\begin{theo}
$\Psi:\,U_{r,s}(C_n^{(1)}) \longrightarrow {\mathcal
U}'_{r,s}(C_n^{(1)})$ is an epimorphism.
\end{theo}

\begin{theo}
{\quad  ${\mathcal
U\,'}_{r,s}(\hat{\frak {g}})={\mathcal U}_{r,s}(\hat{\frak {g}})$.}
\end{theo}

\begin{theo}
There exists a surjective $\Phi:\,{\mathcal U\,'}_{r,s}(\hat{\frak {g}}) \longrightarrow U_{r,s}(\hat{\frak {g}})$ such that $\Psi\Phi=\Phi\Psi=1$.
\end{theo}
Therefore, to prove the Drinfeld Isomorphism Theorem is equivalent
to prove the above three Theorems. We only prove the Theorem 3.13,
and the last two theorems are the same as Theorems B and C in
\cite{HZ}, which are left to the reader.

\medskip

\noindent {\it Proof of Theorem 3.12.}  We shall check that elements
$E_i,~F_i,~\om_i,~\om'_i\; (i\in I_0),\, \gamma^{\pm\frac{1}2}$, $
\gamma'^{\,\pm\frac{1}2},\, D, D'$ satisfy the defining relations of
$(\hat{\rm C}1)-(\hat{\rm C}7)$ of $U_{r,s}(C_n^{(1)})$. At first,
the defining relations of ${\mathcal U}_{r,s}(C_n^{(1)})$ imply that
$E_i,\,F_i,\,\om_i,\,\om'_i \; (i\in I)$ generate a subalgebra of
${\mathcal U}_{r,s}(C_n^{(1)})$ that is isomorphic to
$U_{r,s}({C}_n)$. So we are left to check the relations involving
$i=0$.

For the proof of the relations of
$(\hat{\rm C}1)-(\hat{\rm C}3)$, it is almost the same as those for the simply-laced cases, please see (\cite{HZ}).
Here we only give the proof of the relations of $(\hat{\rm C}4)-(\hat{\rm C}7)$.

For $(\hat{\rm C}4)$: \ at first, when $i\neq 0$, we see that
$$
[\,E_0,F_i\,]
=[\,a\,x_{\theta}^-(1) \,(\gamma'^{-1}\,\om_{\theta}^{-1}),\,
x_i^-{(0)}\,] =- a\,[\,x_i^-{(0)},\,x_{\theta}^-(1) \,]_{\langle
\omega'_i,~ \omega_0 \rangle^{-1}}
(\gamma'^{-1}\,\om_{\theta}^{-1}).
$$

By the $(r,s)$-Serre relations, we claim
\begin{lemm} \ $[\,x_i^-{(0)},\, x_{\theta}^-(1)\,]_{\langle \omega'_i,~ \omega_0
\rangle^{-1}}=0$, for $i\in I$.
\end{lemm}

%For convention, we denote $x_{\alpha_{i_1,\cdots,
%i_n}}^{\pm}(k)=x_{i_1\cdots i_n}^{\pm}(k)$.

\begin{proof} \ (I) \ When $i=1$,
$\langle\om_1',\om_0\rangle=s$ and $\la
\om_1',\om_\theta\ra=s^{-1}$: we first notice that
\begin{equation*}
\begin{split}
x_{\theta}^-(1)=&x_{\beta_{1\,1}}^-(1)=[\,x_1^-(0),\,x_2^-(0)
\,x_{\beta_{1\,3}}^-(1)\,]_{(r^{-\frac{1}{2}},\,(rs)^{-\frac{1}{2}})}
\quad {\hbox{(using (\ref{b:3}))}}\\
=&[\,[\,x_1^-(0),\,x_{2}^-(0)\,]_{r^{-\frac{1}{2}}}
,\,x_{\beta_{1\,3}}^-(1)\,]_{(rs)^{-\frac{1}{2}}}\\
&+r^{-\frac{1}{2}}[\,x_2^-(0),\,
\underbrace{[\,x_{1}^-(0),\,x_{\beta_{1\,3}}^-(1)\,]_{s^{-\frac{1}{2}}}}
\,]\quad{\hbox{(=0 by (\ref{b:18}))}}\\
=&[\,[\,x_1^-(0),\,x_{2}^-(0)\,]_{r^{-\frac{1}{2}}}
,\,x_{\beta_{1\,3}}^-(1)\,]_{(rs)^{-\frac{1}{2}}}.
\end{split}
\end{equation*}
As a result, it is no difficult to see that
\begin{eqnarray*}
\begin{split}
&[\,x_1^-(0),\,x_{\theta}^-(1)\,]_{s^{-1}}\\
&\quad=[\,x_1^-(0),\,[\,x_{1}^-(0),\,x_2^-(0)\,]_{r^{-\frac{1}{2}}}
,\,x_{\beta_{1\,3}}^-(1)\,]_{((rs)^{-\frac{1}{2}},\,s^{-1})}\quad{\hbox{(using (\ref{b:3}))}}\\
&\quad=[\,\underbrace{[\,x_1^-(0),\,[\,x_{1}^-(0),\,x_2^-(0)\,]_{r^{-\frac{1}{2}}}
\,]_{s^{-\frac{1}{2}}}}
,\,x_{\beta_{1\,3}}^-(1)\,]_{r^{-\frac{1}{2}}s^{-1}}\quad{\hbox{(=0 by (\ref{b:10}))}}\\
&\quad+[\,[\,x_{1}^-(0),\,x_2^-(0)\,]_{r^{-\frac{1}{2}}}
,\,\underbrace{[\,x_1^-(0),\,x_{\beta_{1\,3}}^-(1)\,]_{s^{-\frac{1}{2}}}}
\,]_{r^{-\frac{1}{2}}}\quad{\hbox{(=0 by (\ref{b:18}))}}\\
&\quad=0.
\end{split}
\end{eqnarray*}

\ (II) \ When $2\leqslant i \leqslant n-1$,
$\langle\om_{i}',\om_0\rangle=1$, that is, $\la
\om_{i}',\om_\theta\ra=1$. Using (\ref{b:3}) and the $(r,s)$-Serre relations, one
has
\begin{eqnarray*}
\begin{split}
&\quad[\,x_i^-{(0)},\, x_{\theta}^-(1)\,]
\qquad {\hbox{(by definition)}}\\
&\quad=[\,x_i^-(0),\,[\,x_{1}^-(0),\,\cdots, x_{i-2}^-(0),\,
x_{\beta_{1\,i-1}}^-(1)\,]_{(r^{-\frac{1}{2}},\cdots,r^{-\frac{1}{2}}
,\,(rs)^{-\frac{1}{2}})}\,]\\
&\quad=[\,x_{1}^-(0),\,\cdots, x_{i-2}^-(0),\,
\underbrace{[\,x_i^-(0),\,x_{\beta_{1\,i-1}}^-(1)\,]}
\,]_{(r^{-\frac{1}{2}},\cdots,r^{-\frac{1}{2}},\,(rs)^{-\frac{1}{2}})}.
\end{split}
\end{eqnarray*}
Thus, it remains to check that
$[\,x_i^-(0),\,x_{\beta_{1,i-1}}^-(1)\,]=0$.

Actually, it is easy to get
\begin{eqnarray*}
\begin{split}
&\quad[\,x_i^-{(0)},\, x_{\beta_{1\,i-1}}^-(1)\,]_{r^{-\frac{1}{2}}s^{\frac{1}{2}}}
\qquad {\hbox{(by definition)}}\\
&\quad=[\,x_i^-(0),\,[\,x_{i-1}^-(0),\,x_i^-(0),\,x_{\beta_{1\,i+1}}^-(1)
\,]_{(r^{-\frac{1}{2}},\,r^{-\frac{1}{2}})}\,]_{r^{-\frac{1}{2}}s^{\frac{1}{2}}}
\qquad {\hbox{(using (\ref{b:3}))}}\\
&\quad=[\,x_i^-(0),\,[\,x_{i-1}^-(0),\,x_i^-(0)\,]_{r^{-\frac{1}{2}}}
,\,x_{\beta_{1\,i+1}}^-(1)\,]_{(r^{-\frac{1}{2}},\,r^{-\frac{1}{2}}s^{\frac{1}{2}})}
\qquad {\hbox{(using (\ref{b:3}))}}\\
&\qquad+r^{\frac{1}{2}}[\,x_i^-(0),\,x_i^-(0),\,
\underbrace{[\,x_{i-1}^-(0),\,x_{\beta_{1\,i+1}}^-(1)\,]}
\,]_{(1,\,r^{-\frac{1}{2}}s^{\frac{1}{2}})}\quad {\hbox{(=0 by (\ref{b:19}))}}\\
&\quad=[\,\underbrace{[\,x_i^-(0),\,[\,x_{i-1}^-(0),\,x_i^-(0)\,]_{r^{-\frac{1}{2}}}
\,]_{s^{\frac{1}{2}}}}
,\,x_{\beta_{1\,i+1}}^-(1)\,]_{r^{-1}}\quad {\hbox{(=0 by (\ref{b:9}))}}\\
&\qquad+s^{\frac{1}{2}}[\,[\,x_{i-1}^-(0),\,x_i^-(0)\,]_{r^{-\frac{1}{2}}}
,\,[\,x_{i}^-(0),\,x_{\beta_{1\,i+1}}^-(1)\,]_{r^{-\frac{1}{2}}}
\,]_{(rs)^{-\frac{1}{2}}}\quad {\hbox{(using (\ref{b:4}))}}\\
&\quad=s^{\frac{1}{2}}[\,x_{i-1}^-(0),\,
\underbrace{[\,x_i^-(0),\,x_{\beta_{1\,i}}^-(1)\,]_{s^{-\frac{1}{2}}}}
\,]_{r^{-1}}\quad {\hbox{(=0 by (\ref{b:20}))}}\\
&\qquad+[\,[\,x_{i-1}^-(0),\,x_{\beta_{1\,i}}^-(1)\,]_{r^{-\frac{1}{2}}}
,\,x_{i}^-(0)\,]_{r^{-\frac{1}{2}}s^{\frac{1}{2}}}
\quad {\hbox{(by definition)}}\\
&\quad=[\,x_{\beta_{1\,i-1}}^-(1),\,x_{i}^-(0)\,]_{r^{-\frac{1}{2}}s^{\frac{1}{2}}},
\end{split}
\end{eqnarray*}
which implies that
$(1+r^{-\frac{1}{2}}s^{\frac{1}{2}})[\,x_i^-(0),\,
x_{\beta_{1\,i-1}}^-(1)\,]=0.$  Then if $r\ne -s$, we get the
required relation
\begin{equation}
[\,x_i^-(0),\,x_{\beta_{1,i-1}}^-(1)\,]=0.
\end{equation}

\ (III) \, When $i=n$, $\langle\om_{n}',\om_0\rangle=(rs)^{-1}$,
that is, $\la \om_{n}',\om_\theta\ra=rs$. Applying (\ref{b:3}) and the
$(r,s)$-Serre relations, we deduce from direct calculations
\begin{eqnarray*}
\begin{split}
&[\,x_n^-{(0)},\, x_{\theta}^-(1)\,]_{rs}
\qquad {\hbox{(by definition)}}\\
&\quad=[\,x_n^-(0),\,[\,x_{1}^-(0),\,\cdots, x_{n-2}^-(0),\,
x_{\beta_{1\,n-1}}^-(1)\,]_{(r^{-\frac{1}{2}},\cdots,r^{-\frac{1}{2}}
,\,(rs)^{-\frac{1}{2}})}
\,]_{rs}\\
&\quad=[\,x_{1}^-(0),\,\cdots, x_{n-2}^-(0),\,
\underbrace{[\,x_n^-(0),\,x_{\beta_{1\,n-1}}^-(1)\,]_{rs}}
\,]_{(r^{-\frac{1}{2}},\cdots,r^{-\frac{1}{2}},\,(rs)^{-\frac{1}{2}})}.
\end{split}
\end{eqnarray*}
It remains  to show the relation
$[\,x_n^-(0),\,x_{\beta_{1\,n-1}}^-(1)\,]_{rs}=0$.

In fact,
\begin{eqnarray*}
\begin{split}
&[\,x_n^-{(0)},\, x_{\beta_{1\,n-1}}^-(1)\,]_{s^{2}}
\qquad {\hbox{(by definition)}}\\
&\quad=[\,x_n^-(0),\,[\,x_{n-1}^-(0),\,x_n^-(0),\,x_{\alpha_{1\,n-1}}^-(1)
\,]_{(s,\,r^{-\frac{1}{2}})}\,]_{s^{2}}
\qquad {\hbox{(using (\ref{b:7}))}}\\
&\quad=[\,x_n^-(0),\,[\,x_{n-1}^-(0),\,x_n^-(0)\,]_{r^{-1}}
,\,x_{\alpha_{1\,n-1}}^-(1)\,]_{(r^{\frac{1}{2}}s,\,s^{2})}
\qquad {\hbox{(using (\ref{b:7}))}}\\
&\quad+r^{-1}[\,x_n^-(0),\,x_n^-(0),\,
\underbrace{[\,x_{n-1}^-(0),\,x_{\alpha_{1\,n-1}}^-(1)\,]_{r^{\frac{1}{2}}}}
\,]_{(rs,\,s^{2})}\quad {\hbox{(=0 by (\ref{b:16}))}}\\
&\quad=[\,\underbrace{[\,x_n^-(0),\,[\,x_{n-1}^-(0),\,x_n^-(0)\,]_{r^{-1}}
\,]_{s}}
,\,x_{\alpha_{1\,n-1}}^-(1)\,]_{r^{\frac{1}{2}}s^2}\quad {\hbox{(=0 by (D9$_2$))}}\\
&\quad+s[\,[\,x_{n-1}^-(0),\,x_n^-(0)\,]_{r^{-1}}
,\,[\,x_{n}^-(0),\,x_{\alpha_{1\,n-1}}^-(1)\,]_{s}
\,]_{r^{\frac{1}{2}}}\qquad {\hbox{(by definition \& (\ref{b:4}))}}\\
&\quad=s[\,x_{n-1}^-(0),\,
\underbrace{[\,x_n^-(0),\,x_{\alpha_{1\,n}}^-(1)\,]_{r}}
\,]_{r^{-\frac{3}{2}}}\quad {\hbox{(=0 by (\ref{b:21}))}}\\
&\quad+rs[\,[\,x_{n-1}^-(0),\,x_{\alpha_{1\,n}}^-(1)\,]_{r^{-\frac{1}{2}}}
,\,x_{n}^-(0)\,]_{r^{-2}}
\quad {\hbox{(by definition)}}\\
&\quad=rs[\,x_{\beta_{1\,n-1}}^-(1),\,x_{n}^-(0)\,]_{r^{-2}}.
\end{split}
\end{eqnarray*}
Expanding the two sides of the above result, one has
$$(1+r^{-2}s^{2})[\,x_n^-(0),\, x_{\beta_{1\,n-1}}^-(1)\,]_{rs}=0.$$
Therefore, if $r\ne -s$,
$[\,x_n^-(0),\,x_{\beta_{1\,n-1}}^-(1)\,]_{rs}=0.$

This complete the proof of Lemma 3.15.
\end{proof}

Next, we turn to check the relation
\begin{lemm} \
$[\,E_0,\, F_0\,]=\frac{\omega_0-\omega'_0}{r-s}$.
\end{lemm}
\begin{proof}\,
Using $(\textrm{D1})$ and $(\textrm{D5})$, one has
\begin{equation}
\begin{split}
\bigl[\,E_0,
F_0\,\bigr]&=(rs)^{\frac{n-2}{2}}[2]_1^{-1}\bigl[\,x^-_{\theta}(1)\,
\gamma'^{-1}{\om_\theta}^{-1},\,\gamma^{-1}{\om'_\theta}^{-1}
x^+_{\theta}(-1)\,\bigr]\\
&=(rs)^{\frac{n-2}{2}}[2]_1^{-1}\bigl[\,x^-_{\theta}(1),\,x^+_{\theta}(-1)\,\bigr]\,
\cdot(\gamma^{-1}\gamma'^{-1}{\om_\theta}^{-1}{\om'_\theta}^{-1}).
\end{split}
\end{equation}
First recall the notations,
\begin{gather*}
x^-_{\theta}(1)=[\,x_1^-(0),\cdots, x_n^-(0),\,x_{\alpha_{1\,n-1}}^-(1)\,]
_{(s,\,r^{-\frac{1}{2}},\,\cdots,\,r^{-\frac{1}{2}},\,(rs)^{-\frac{1}{2}})},\\
x_{\theta}^+(-1)=[\,[\,[\,[\,x_{\alpha_{1\,n-1}}^+(-1),\,x_n^+(0)\,]_r,\,
x_{n-1}^+(0)\,]_{r^{-\frac{1}{2}}},\cdots, x_2^+(0)\,]_{r^{-\frac{1}{2}}}
,\,x_1^+(0)\,]_{(rs)^{-\frac{1}{2}}}.
\end{gather*}
Owing to the result of the case of $A_{n-1}^{(1)}$ (\cite{HRZ}), one has
\begin{equation}[\,x_{\alpha_{1\,n-1}}^-(1),\,x_{\alpha_{1\,n-1}}^+(-1)\,]=
\frac{\gamma\om_{\alpha_{1\,n-1}}'-\gamma'\om_{\alpha_{1\,n-1}}}
{r^{\frac{1}{2}}-s^{\frac{1}{2}}}.
\end{equation}
Next, it is easy to get
\begin{equation*}
\begin{split}
&[\,x_{\alpha_{1\,n}}^-(1),\,x_{\alpha_{1\,n}}^+(-1)\,]\qquad{\hbox{(by definition)}}\\
&\quad=[\,[\,x_{n}^-(0),\,x_{\alpha_{1\,n-1}}^-(1)\,]_{s}
,\,[\,x_{\alpha_{1\,n-1}}^+(-1),\,x_{n}^+(0)\,]_{r}\,]\qquad{\hbox{(using (\ref{b:3}))}}\\
&\quad=[\,[\,[\,x_{n}^-(0),\,x_{\alpha_{1\,n-1}}^+(-1)\,],\,x_{\alpha_{1\,n-1}}^-(1)\,]_{s}
,\,x_{n}^+(0)\,]_{r}\quad{\hbox{(=0 by (\ref{b:3}) \& (D8))}}\\
&\quad+[\,[\,x_{n}^-(0),\,[\,x_{\alpha_{1\,n-1}}^-(1),\,x_{\alpha_{1\,n-1}}^+(-1)\,]\,]_{s}
,\,x_{n}^+(0)\,]_{r}\quad{\hbox{(using (3.24), (D5), (D8))}}\\
&\quad+[\,x_{\alpha_{1\,n-1}}^+(-1),\,[\,
[\,x_{n}^-(0),\,x_{n}^+(0)\,],\,x_{\alpha_{1\,n-1}}^-(1)\,]_{s}\,]_{r}\quad{\hbox{(using (D8), (D5), (3.24))}}\\
&\quad+[\,x_{\alpha_{1\,n-1}}^+(-1),\,[\,x_{n}^-(0),\,
[\,x_{\alpha_{1\,n-1}}^-(1),\,x_{n}^+(0)\,]\,]_{s}
\,]_{r}\quad{\hbox{(=0 by (\ref{b:3}) \& (D8))}}\\
&\quad=\gamma\om'_{\alpha_{1\,n-1}}\cdot \frac{\om'_n-\om_n}
{r^{\frac{1}{2}}-s^{\frac{1}{2}}}
+\frac{\gamma\om'_{\alpha_{1\,n-1}}-\gamma'\om_{\alpha_{1\,n-1}}}
{r^{\frac{1}{2}}-s^{\frac{1}{2}}}\om_n\\
&\quad=\frac{\gamma\om'_{\alpha_{1\,n}}-\gamma'\om_{\alpha_{1\,n}}}
{r^{\frac{1}{2}}-s^{\frac{1}{2}}}.
\end{split}
\end{equation*}
Furthermore, we obtain
\begin{equation*}
\begin{split}
&[\,x_{\beta_{1\,n-1}}^-(1),\,x_{\beta_{1\,n-1}}^+(-1)\,]\qquad{\hbox{(by definition)}}\\
&\quad=[\,[\,x_{n-1}^-(0),\,x_{\alpha_{1\,n}}^-(1)\,]_{r^{-\frac{1}{2}}}
,\,[\,x_{\alpha_{1\,n}}^+(-1),\,x_{n-1}^+(0)\,]_{s^{-\frac{1}{2}}}\,]
\qquad{\hbox{(using (\ref{b:3}))}}\\
&\quad=[\,[\,[\,x_{n-1}^-(0),\,x_{\alpha_{1\,n}}^+(-1)\,]
,\,x_{\alpha_{1,n}}^-(1)\,]_{r^{-\frac{1}{2}}}
,\,x_{n-1}^+(0)\,]_{s^{-\frac{1}{2}}}\quad{\hbox{(=0 by (\ref{b:3}) \& (D5))}}\\
&\quad+[\,[\,x_{n-1}^-(0),\,
[\,x_{\alpha_{1\,n}}^-(1),\,x_{\alpha_{1\,n}}^+(-1)\,]\,]_{r^{-\frac{1}{2}}}
,\,x_{n-1}^+(0)\,]_{s^{-\frac{1}{2}}}\quad{\hbox{(using (D5) \& (D8))}}\\
&\quad+[\,x_{\alpha_{1\,n}}^+(-1),\,[\,
[\,x_{n-1}^-(0),\,x_{n-1}^+(0)\,],\,x_{\alpha_{1\,n}}^-(1)\,]_{r^{-\frac{1}{2}}}
\,]_{s^{-\frac{1}{2}}}\quad{\hbox{(using (D8) \& (D5))}}\\
&\quad+[\,x_{\alpha_{1\,n}}^+(-1),\,[\,x_{n-1}^-(0),\,
[\,x_{\alpha_{1\,n}}^-(1),\,x_{n-1}^+(0)\,]\,]_{r^{-\frac{1}{2}}}
\,]_{s^{-\frac{1}{2}}}\quad{\hbox{(=0 by (\ref{b:3}) \& (D5))}}\\
&\quad=(rs)^{-\frac{1}{2}}\gamma\om'_{\alpha_{1\,n}}
\cdot \frac{\om'_{n-1}-\om_{n-1}}{r^{\frac{1}{2}}-s^{\frac{1}{2}}}
+(rs)^{-\frac{1}{2}}\frac{\gamma\om'_{\alpha_{1\,n}}-\gamma'\om_{\alpha_{1\,n}}}
{r^{\frac{1}{2}}-s^{\frac{1}{2}}}\om_{n-1}\\
&\quad=(rs)^{-\frac{1}{2}}\frac{\gamma\om'_{\beta_{1\,n-1}}-\gamma'\om_{\beta_{1\,n-1}}}
{r^{\frac{1}{2}}-s^{\frac{1}{2}}}.
\end{split}
\end{equation*}
Repeating the above step, we get by directly calculating
\begin{equation}
\begin{split}
[\,x_{\beta_{1\,2}}^-(1),\,x_{\beta_{1\,2}}^+(-1)\,]
=(rs)^{\frac{2-n}{2}}\frac{\gamma\om'_{\beta_{1\,2}}-\gamma'\om_{\beta_{1\,2}}}
{r^{\frac{1}{2}}-s^{\frac{1}{2}}}.
\end{split}
\end{equation}
As a result, we get the required conclusion:
\begin{equation*}
\begin{split}
&[\,x_{\beta_{1\,1}}^-(1),\,x_{\beta_{1\,1}}^+(-1)\,]\qquad{\hbox{(by definition)}}\\
&\quad=[\,[\,x_{1}^-(0),\,x_{\beta_{1\,2}}^-(1)\,]_{(rs)^{-\frac{1}{2}}}
,\,[\,x_{\beta_{1\,2}}^+(-1),\,x_{1}^+(0)\,]_{(rs)^{-\frac{1}{2}}}\,]
\qquad{\hbox{(using (\ref{b:3}))}}\\
&\quad=[\,[\,[\,x_{1}^-(0),\,x_{\beta_{1\,2}}^+(-1)\,],\,x_{\beta_{1\,2}}^-(1)
\,]_{(rs)^{-\frac{1}{2}}}
,\,x_{1}^+(0)\,]_{(rs)^{-\frac{1}{2}}}\quad{\hbox{(using (\ref{b:3}), (D8) \& (D5))}}\\
&\quad+[\,[\,x_{1}^-(0),\,[\,x_{\beta_{1\,2}}^-(1),\,x_{\beta_{1\,2}}^+(-1)\,]
\,]_{(rs)^{-\frac{1}{2}}}
,\,x_{1}^+(0)\,]_{(rs)^{-\frac{1}{2}}}\quad{\hbox{(=0 by (D8) \& (D5))}}\\
&\quad+[\,x_{\beta_{1\,2}}^+(-1),\,[\,
[\,x_{1}^-(0),\,x_{1}^+(0)\,],\,x_{\beta_{1\,2}}^-(1)\,]_{(rs)^{-\frac{1}{2}}}
\,]_{(rs)^{-\frac{1}{2}}}\quad{\hbox{(=0 by (D8) \& (D5))}}\\
&\quad+[\,x_{\beta_{1\,2}}^+(-1),\,[\,x_{1}^-(0),\,
[\,x_{\beta_{1\,2}}^-(1),\,x_{1}^+(0)\,]\,]_{(rs)^{-\frac{1}{2}}}
\,]_{(rs)^{-\frac{1}{2}}}\quad{\hbox{(using (\ref{b:3}), (D8) \& (D5))}}\\
&\quad=(rs)^{\frac{2-n}{2}}\gamma\om'_{\beta_{1\,2}}\cdot \frac{\om'_{1}-\om_{1}}
{r^{\frac{1}{2}}-s^{\frac{1}{2}}}
+(rs)^{\frac{2-n}{2}}\frac{\gamma\om'_{\beta_{1\,2}}-\gamma'\om_{\beta_{1\,2}}}
{r^{\frac{1}{2}}-s^{\frac{1}{2}}}\om_{1}\\
&\quad=(rs)^{\frac{2-n}{2}}\frac{\gamma\om'_{\beta_{1\,1}}-\gamma'\om_{\beta_{1\,1}}}
{r^{\frac{1}{2}}-s^{\frac{1}{2}}}.
\end{split}
\end{equation*}
Thus, we arrive at the last step
$$\bigl[\,E_0,
F_0\,\bigr]=\frac{\gamma'^{-1}\om_{\theta}^{-1}-\gamma^{-1}\om_{\theta'}^{-1}}{r-s}.$$

The proof is complete.
\end{proof}

For $(\hat{\rm C}6)$: \ when $i\cdot j \neq 0$, $(D9)$ implies that
the corresponding generators satisfy exactly those $(r,s)$-Serre
relations in $U_{r,s}(\mathrm{C}_n^{(1)})$, so it is enough to check
the $(r,s)$-Serre relations involving $i\cdot j=0$.

\begin{lemm} {\it $(1)$ \ $E_nE_0=rs\,E_0E_n$,

$(2)$ \ $E_0^2E_1-(r+s)E_0E_1E_0+rsE_1E_0^2=0$,

$(3)$ \ $E_0E_1^3-(r+(rs)^{\frac{1}{2}}+s)E_1E_0E_1^2+
(rs)^{\frac{1}{2}}(r+(rs)^{\frac{1}{2}}+s)E_1^2E_0E_1
-(rs)^{\frac{3}{2}}E_1^3E_0=0$,

$(4)$ \ $F_0F_n=rs\,F_nF_0$,

$(5)$ \ $F_{1}F_0^{2}-(r+s)\,F_{0}F_1F_{0}+rs\,F_0^2F_{1}=0$,

$(6)$ \ $F_1^3F_{0}-(r+(rs)^{\frac{1}{2}}+s)
F_1^2F_{0}F_1+(rs)^{\frac{1}{2}}(r+(rs)^{\frac{1}{2}}+s)F_1F_0F_{1}^2
-(rs)^{\frac{3}{2}}F_0F_1^3=0$.}
\end{lemm}

\begin{proof} \,  Here we only check the first and third $(r,s)$-Serre
relations, and the rest are left to the readers.

(1) \,  We start with the following calculations
\begin{eqnarray*}
\begin{split}
&E_nE_0-rsE_0E_n\qquad \hbox{(by definition)}\\
&\quad=[\,x_n^+(0),\,x_{\theta}^-(1)\gamma'^{-1}\om_{\theta}^{-1}\,]_{rs}
\qquad \hbox{(using (D5))}\\
&\quad=[\,x_n^+(0),\,x_{\theta}^-(1)\,]\cdot \gamma'^{-1}\om_{\theta}^{-1}
\qquad \hbox{(by definition)}\\
&\quad=[\,x_1^-(0),\,\cdots,\,x_{n-1}^-(0),\underbrace{[\,x_n^+(0),\,x_n^-(0)\,]},\,x_{\alpha_{1\,n-1}}^-(1)\,]
_{(s,\,r^{-\frac{1}{2}},\,\cdots,\,r^{-\frac{1}{2}},\,(rs)^{-\frac{1}{2}})}\\
&\quad=[\,x_1^-(0),\,\cdots,\,\underbrace{[\,x_{n-1}^-(0),\,x_{\alpha_{1\,n-1}}^-(1)\,]_{r^{\frac{1}{2}}}}\,]
_{(r^{-\frac{1}{2}},\,\cdots,\,r^{-\frac{1}{2}},\,(rs)^{-\frac{1}{2}})}\om_n
\quad \hbox{(=0 by (\ref{b:16}))}\\
&\quad=0.
\end{split}
\end{eqnarray*}

(2) \,  At the same time, we consider that
\begin{equation*}
\begin{split}
&\bigl[\,E_1, x^-_{\theta}(1)\,\bigr]\qquad \hbox{(using (\ref{b:3}))}\\
&\quad=[\,[\,x_1^{+}(0),\,x_1^-(0)\,],\,x_{\beta_{1\,2}}^-(1)\,]_{(rs)^{-\frac{1}{2}}}
\qquad(\hbox{=0 by (D8) \& (D5)})\\
&\quad+[x_1^-(0),\cdots,x_n^-(0),\cdots,x_2^-(0),
[\,x_1^+(0),x_1^-(1)\,]\,]_{(s^{\frac{1}{2}},\cdots,s^{\frac{1}{2}},s,r^{-\frac{1}{2}}
\cdots,r^{-\frac{1}{2}},\,(rs)^{-\frac{1}{2}})}\\
&\quad=-(rs)^{\frac1{4}}(\gamma\gamma')^{-\frac{1}{2}}\cdot
[\,x_1^-(0),x_{\beta_{2\,2}}^-(1)\,]_{r^{-1}}\,\om_1.
\end{split}
\end{equation*}
In terms of the above result, it is easy to see that
\begin{equation*}
\begin{split}
 &E_0E_1^3-(r+(rs)^{\frac{1}{2}}+s)E_1E_0E_1^2+
(rs)^{\frac{1}{2}}(r+(rs)^{\frac{1}{2}}+s)E_1^2E_0E_1-(rs)^{\frac{3}{2}}E_1^3E_0\\
&\quad=(rs)^{\frac{3}{2}}\Big(E_1^3 x^-_{\theta}(1)-(1{+}(r^{-1}s)^{\frac{1}{2}}{+}
r^{-1}s)E_1^2x^-_{\theta}(1)E_1\\
&\quad\hskip0.5cm +((r^{-1}s)^{\frac{1}{2}}{+}
r^{-1}s{+}(r^{-1}s)^{\frac{3}{2}})E_1x^-_{\theta}(1)E_1^2
-(r^{-1}s)^{\frac{3}{2}}\,x^-_{\theta}(1)E_1^3\Big)
(\gamma'^{-1}\om_{\theta}^{-1})\\
&\quad= (rs)^{\frac{3}{2}}\,\bigl[\,E_1, \,E_1,\,
\underbrace{\bigl[\,E_1,\,x^-_{\theta}(1)\,
\bigr]}\,\bigr]_{((r^{-1}s)^{\frac{1}{2}},\,r^{-1}s)}
\,(\gamma'^{-1}\om_{\theta}^{-1})\\
&\quad=-(rs)^{\frac{3}{4}}[2]_1(\gamma\gamma')^{-\frac{1}{2}}
\underbrace{\bigl[\,E_1,\,x^-_{\beta_{2\,2}}(1)\,
\bigr]}\om_1^2\,(\gamma'^{-1}\om_{\theta}^{-1})\quad
(\hbox{=0 by (\ref{b:3}) \& $(\textrm{D8})$})\\
&\quad=0.
\end{split}
\end{equation*}

The proof is complete.
\end{proof}

For $(\hat{\rm C}7)$: \ the verification is analogous to that of
$(\hat{\rm C}6)$. Hence, we establish the Drinfeld isomorphism for
the two-parameter quantum affine algebra of type $C^{(1)}$.

 \vskip 0.5cm

\section{Vertex representation }
\setcounter{equation}{0}

\subsection{Enlarged quantum Heisenberg algebra}\,  In order to obtain the vertex
representation of $U_{r,s}(C_n^{(1)})$, first we need the
following construction.

The associative algebra $U_{r,s}(\widehat{\frak{h}})$, generated by $\{\,a_i(m),\,
\gamma^{\pm\frac1{2}},\,\gamma'^{\pm\frac1{2}} \mid
m\in\mathbb{Z}\backslash \{0\},\,i$ $=1,\,2\,\cdots, n\}$, is called the two-parameter quantum Heisenberg algebra, where $a_i(m)$ $(i=1,\,2, \cdots, n)$ and $\omega_j^{\pm1}$, $\omega_j'^{\pm 1}$'s satisfy $(D3)$, as well as the following
\begin{eqnarray}
\qquad [a_i(m), a_j(l)]=\delta_{m+l, 0} \frac{
(\gamma\gamma')^{\frac{|m|}{2}}(rs)^{-\frac{m (\alpha_i|\alpha_j)}2}[\,m
(\alpha_i|\alpha_j)\,]}{|m|}\cdot\frac{\gamma^{|m|}-\gamma'^{|m|}}{r-s}.
\end{eqnarray}

Define an enlarged quantum Heisenberg algebra $U_{r,s}(\widetilde {\mathfrak h})$,
which is an associative algebra generated by the subalgebra $U_{r,s}(\widehat{\frak{h}})$, together with
a coupled bosonic operators $b(m)$ ($m\in\mathbb Z^*$) satisfying the following relations:
\begin{equation}
\begin{array}{rcl}
&& [b_i(m), b_j(l)]=[a_i(m), a_j(l)],\\
&&[a_i(m), b_j(l)]=0, \\
&& [a_i(0), a_j(m)]=[b_i(0), b_j(m)]=0.
\end{array}
\end{equation}

\subsection{Fock space}\,  In order to construct the Fock space, one can take a copy
of the root lattice of $A_{n-1}$ as the sublattice
$\tilde{Q}=Q[A_{n-1}]$ of short roots of the root lattice $Q$ for type $C_n$. The basis of
$\tilde{Q}$ will be denoted by $\tilde{\a}_i, i=1, \cdots, n-1$.
Thus
$$(\tilde{a}_i|\tilde{a}_j)=(\alpha_i|\alpha_j)=\delta_{ij}-\frac12\delta_{|i-j|, 1}.
$$
We also consider the associated weight lattice
$\tilde{P}=P[A_{n-1}]$ defined by the inner product.

The Fock space ${\mathcal{V}}$ is defined to be the tensor product
of the symmetric algebra generated by $a_i(-m)$'s, $b_i(-m)$'s $(m\ne 0)$ and the
group algebra generated by $e^{\lambda}\otimes e^{\tilde{\lambda}}$
such that $(a_i|\lambda)\pm(\tilde{a}_i|\tilde{\lambda})\in
{\bf Z}$ for each $i\in\{ 1, \cdots, n\}$, where $\lambda\in P$ and
$\tilde{\lambda}\in \tilde{P}$. Note that we treat $\tilde{a}_n=0$.

The actions of $a_i(m)$'s and $b_i(m)$'s with $m\neq 0$ on
$\mathcal{V}$ are obtained by viewing the Fock space $\mathcal{V}
$ as some quotient space of the Heisenberg algebra tensored with the
group algebras of $P$ and $\tilde{P}$. The operators $a_i(0)$,
$b_i(0)$, $e^{\a}$, $e^{\tilde{a}}$ act on $\mathcal{V} $ by
\begin{eqnarray*}
a_i(0)e^{\lambda}e^{\tilde{\lambda}}&=&(\a_i|\lambda)e^{\lambda}e^{\tilde{\lambda}},\quad
b_i(0)e^{\lambda}e^{\tilde{\lambda}}=(\tilde{a}_i|\tilde{\lambda})e^{\lambda}e^{\tilde{\lambda}},\\
e^{\a}e^{\lambda}e^{\tilde{\lambda}}&=&e^{\a+\lambda}e^{\tilde{\lambda}},\qquad\quad\ \;
e^{\tilde{a}}e^{\lambda}e^{\tilde{\lambda}}=e^{\lambda}e^{\tilde{a}+\tilde{\lambda}};\\
\om_i\cdot(v\otimes e^\be)&=&\la \be,\, i\ra \,v\otimes e^\be,\qquad
\om_i'\cdot (v\otimes
e^\be)=\la i,\, \be\ra^{-1}\,v\otimes e^\be
\end{eqnarray*}
and define
$$D(r)(v\otimes e^\be)=r^{\be}v\otimes e^\be,\quad D(s)(v\otimes e^\be)=s^{\be}v\otimes e^\be.$$

It is easy to see that $a_i(m), b_j(l), e^{\a_i}, e^{\tilde{a}_j}$
commute with each other except that
$$
[a_i(0), e^{\alpha_j}]=(\alpha_i|\alpha_j)e^{\alpha_j} , \qquad [b_i(0),
e^{\tilde{a}_j}]=(\tilde{a}_i|\tilde{a}_j)e^{\tilde{a}_j}.
$$
\subsection{Normal ordering}\, The normal product $: \quad :$ is defined as usual:
   $$ :a_i(m) a_j(l): = a_i(m) a_j(l) \quad (\textit{if } m \leq l),
   \; \mbox{ or } \ a_j(l) a_i(m) \quad (\textit{if } m>l), $$
   $$ :e^{\a} a_i(0):=:a_i(0) e^{\a}:=e^{\a} a_i(0) \, , $$
   $$ :e^{\tilde{a}\tilde{a}} b_i(0):=:b_i(0) e^{\tilde{a}}:=e^{\tilde{a}} b_i(0) \, , $$
   and similarly for product involving the $b_i(m)$.

\subsection{Quasi-cocycle}\, Let $\vep( \ \ ,\ \ )$: $P\times P\longrightarrow \mathbb{K}$
 be the quasi-cocycle such that
\begin{eqnarray*}
\vep(\a, \be+\theta)&=&\vep(\a,\be)\vep(\a, \theta),\\
\vep(\a+\be, \theta)&=&\vep(\a, \theta)\vep(\be, \theta)
(-1)^{(\overline{\a+\be}-\overline{\a}-\overline{\be}|\overline{\theta})}.
\end{eqnarray*}
where the $-$ is the projection from $P$ to $\tilde P$ defined by
$$\a=\sum_{i=1}^n m_i\a_i\in P \longmapsto \overline{\a_i}=
\sum_{i=1}^{n-1} \overline{m_i}\a_i, \quad \overline{m_i}=m_i\ (\mod\,
2).
$$
We construct such a quasi-cocycle directly by
$$\epsilon(\alpha_i,\alpha_j)=\left\{\begin{array}{cl} (-1)^{a_{ij}}(r_is_i)^{\frac{a_{ij}}{2}},
&\alpha_i+\alpha_j\in\Phi, \ i> j;\\
-(r_is_i)^{\frac{1}{2}},
& \ i= j;\\
1,& ~\hbox{other pairs}~ (i,j).
\end{array}\right.$$

It is easy to verify that the quasi-cocycle satisfies all the
defining relations. In particular, we have
\begin{equation}
\vep(\a_i,\a_j)\vep(\a_j,\a_i)=\left\{\begin{array}{ll}
                (-1)^{2(\a_i|\a_j)}(r_is_i)^{\frac{a_{ij}}{2}}, & \mbox{if}\ \ 1\leq i, j\leq n-1;\\
                -(rs)^{\frac{1}{2}},  & \mbox{otherwise}.
                \end{array}\right.
\end{equation}

For $\a\in P$, we define the operators $\vep_{\a}$ on $\mathcal{V}$
such that
\begin{equation}
\vep_{\a}e^{\lambda}e^{\tilde{\lambda}}=\vep(\a,
\lambda)e^{\lambda}e^{\tilde{\lambda}}.
\end{equation}
For simplicity, we denote $\vep_i=\vep_{\a_i}$ for $i=1, \cdots, n$.

\subsection{Vertex operators}\, We can now introduce the main vertex operators:
\begin{eqnarray*}
&Y_i^{\pm}(z)&=
         \exp (\pm  \sum^{\infty}_{k=1}
                 \frac{a_i(-k)}{[k]} s^{\pm\frac{k}{2}} z^k)\\
    &&\qquad     \times\exp ( \mp \sum^{\infty}_{k=1}
           \frac{a_i(k)   }{[k]} r^{\mp \frac{k}{2}} z^{-k})
         e^{\pm\a_i} z^{\pm a_i(0)}\vep_i ,\\
%&Y_i^{-}(z)&=
  %       \exp ( -\sum^{\infty}_{k=1}
               %  \frac{a_i(-k)}{[k]} s^{-\frac{k}{2}} z^k)\\
    %&&\qquad     \times\exp ( +\sum^{\infty}_{k=1}
       %    \frac{a_i(k)   }{[k]} r^{\frac{k}{2}} z^{-k})
      %   e^{-\a_i} z^{-a_i(0)}\vep_i ,\\
&U_j^+(z)&=
         \exp ( \sum^{\infty}_{k=1}
                 \frac{b_j(-k)}{[k]}z^k)
         \exp ( - \sum^{\infty}_{k=1}
           \frac{b_j(k)   }{[k]}z^{-k})
         e^{ \tilde{a}_j} z^{ b_j(0)}(rs)^{-\frac{1}{8}},\\
&U_j^-(z)&=
         \exp (-\sum^{\infty}_{k=1}
                 \frac{b_j(-k)}{[k]}z^k)
         \exp (\sum^{\infty}_{k=1}
           \frac{b_j(k)   }{[k]}  z^{-k})
         e^{-\tilde{a}_j} z^{- b_j(0)}(rs)^{-\frac{1}{8}},\\
&Z^{+}_j(z)&=U_j^+(s^{1/2}z)+(-1)^{2a_j(0)}U_j^-(r^{1/2}z),\\
&Z^{-}_j(z)&=U_j^+(r^{-1/2}z)+(-1)^{2a_j(0)}U_j^-(s^{-1/2}z),
\end{eqnarray*}
where $i\in\{1, \cdots, n\}, j\in\{1, \cdots, n-1\}$. For simplicity,
we define $Z_n^{\pm}(z)=1$.

\subsection{Vertex representation of $U_{r,s}(\mathrm{C}_n^{(1)})$}

\begin{theo} \label{T:1} The Fock space
    $\mathcal{V}$ is a $U_{r,s}(\mathrm{C}_n^{(1)})$-module of level $1$
    under the action defined by
\begin{eqnarray*}
\gamma^{\pm\frac1{2}}&\mapsto& r^{\pm\frac1{2}},\qquad \gamma'^{\pm\frac1{2}}\mapsto s^{\pm\frac1{2}},\\
D &\mapsto& D(r),\qquad
D'\mapsto D(s),\\
\omega_i&\mapsto & \omega_i,\qquad\quad\,\,
\omega'_i\mapsto  \omega'_i,\\
 a_{i}(m) &\longmapsto& a_i(m),\\
x_i^{\pm}(z) &\longmapsto&
             Y_i^{\pm}(z)Z_i^{\pm}(z),  \qquad i=1, \cdots, n.
\end{eqnarray*}
\end{theo}

\noindent{\bf Proof of Theorem 4.1.} \ In what follows,
we check Theorem 4.1 using the similar but more detailed techniques than in the one-parameter setting (cf. \cite{JKM2}).
Let $X_i^{\pm}(z)$ denote the images of $x_i^{\pm}(z)$ ($i\in \{1,\,2,\cdots, n\}$)
in the algebra ${\mathcal U}_{r,s}(C_n^{(1)})$ under the mapping $\pi$, respectively.
Therefore we have to check $X_i^{\pm}(z)$, $\Phi_j(z)$ and $\Psi_j(z)$ satisfy
relations $({\rm D}1)$---$({\rm D}9)$. Firstly, it is clear that relations $({\rm D}1)$---$({\rm D}6)$
are true by the construction, which can be verified similarly to \cite{HZ}. We are left to show relations $({\rm D}7)$---$({\rm D}9)$.

Firstly, we give the notation formally, for $a \in \mathbb{K}$,
\begin{eqnarray}{\label{n:1}
(1-z)_{r,s}^{a}:=exp(-\sum_{{n\geq 1}}\frac{[an]}{n[n]}z^{n}).}
\end{eqnarray}
In particular, \begin{eqnarray}
&(1-z)_{r,s}^{1}=(1-z),\\
&(1-z)_{r,s}^{-1}=\Big(\big(1-(rs)^{-1}z\big)\Big)^{-1}.
\end{eqnarray}

\begin{lemm} {\label{l:1} With the above notations, one has
\begin{eqnarray}
&\big(1-(rs)^{\frac{3}{4}}z\big)_{r,s}^{-\frac{1}{2}}\cdot \big(1-(rs)^{\frac{3}{4}}(r^{-1}s)^{\frac{1}{2}}z\big)_{r,s}^{-\frac{1}{2}}
=\big(1-(r^{-1}s)^{\frac{1}{4}}z\big)^{-1};\\
&\big(1-(rs)^{\frac{1}{4}}(rs^{-1})^{\frac{1}{2}}z\big)_{r,s}^{\frac{1}{2}}\big(1-(rs)^{\frac{3}{4}}(r^{-1}s)^{\frac{1}{2}}z\big)_{r,s}^{-\frac{1}{2}}
=\frac{1-(rs^{-1})^\frac{1}{4}z}{1-(r^{-1}s)^\frac{1}{4}z};\\
&\big(1-(rs)^{\frac{1}{4}}a z\big)_{r,s}^{\frac{1}{2}}\cdot \big(1-(rs)^{\frac{3}{4}}a z\big)_{r,s}^{-\frac{1}{2}}
=1;\\
&\big(1-(rs)^{\frac{1}{4}}(rs)^{-\frac{1}{2}}z\big)_{r,s}^{\frac{1}{2}}\big(1-(rs)^{\frac{1}{4}}s^{-1}z\big)_{r,s}^{\frac{1}{2}}
=\big(1-(rs)^{-\frac{1}{4}}s^{-\frac{1}{2}}z\big);\\
&\big(1-(rs)^{\frac{1}{4}}(rs)^{-\frac{1}{2}}z\big)_{r,s}^{\frac{1}{2}}\big(1-(rs)^{\frac{1}{4}}r^{-1}z\big)_{r,s}^{\frac{1}{2}}
=\big(1-(rs)^{-\frac{1}{4}}r^{-\frac{1}{2}}z\big).
\end{eqnarray}
}\end{lemm}
\begin{proof}\, (4.8): using the notation (\ref{n:1}), one directly gets
\begin{eqnarray*}
\begin{split}
LHS&=exp\big(-\sum_{{n\geq 1}}\frac{[-\frac{1}{2}n]}{n[n]}(rs)^{\frac{3n}{4}}z^{n}\big)\cdot
exp\big(-\sum_{{n\geq 1}}\frac{[-\frac{1}{2}n]}{n[n]}(rs)^{\frac{3n}{4}}(r^{-1}s)^{\frac{n}{2}}z^{n}\big)\\
&=exp\big(-\sum_{{n\geq 1}}\frac{[-\frac{1}{2}n]}{n[n]}(rs)^{\frac{3n}{4}}(1+(r^{-1}s)^{\frac{n}{2}})z^{n}\big)\\
&=exp\big(\sum_{{n\geq 1}}\frac{r^{-\frac{n}{4}}s^{\frac{n}{4}}}{n}z^{n}\big)=RHS.
\end{split}
\end{eqnarray*}
The others are similar, which are left to the readers.
\end{proof}
For the later use, we will write
$$\big(z-w\big)_{r,s}^{a}=
\big(1-\frac{w}{z}\big)_{r,s}^{a}\cdot z^a.$$

Before checking relations $({\rm D}7)$--$({\rm D}9)$, we list the operator product expansions.
\begin{lemm} {\label{l:2} With the above notation, the following formulas hold
\begin{eqnarray}
&&Y_i^{\pm}(z)Y_j^{\pm}(w)=:Y_i^{\pm}(z)Y_j^{\pm}(w): \nonumber \\
&&\hskip2.7cm \cdot \left\{\begin{array}{lr}1, &\ (\alpha_i|\alpha_j)=0;\vspace{1mm}\\
\big(z-(rs)^{\frac{3}{4}}(r^{-1}s)^{\pm\frac{1}{2}}w\big)_{r,s}^{-\frac{1}{2}}\cdot\epsilon(\alpha_i,\alpha_j)^{\pm1},
& \ (\alpha_i|\alpha_j)=-\frac{1}{2};\vspace{1mm}\\
\big(z-(r^{-1}s)^{\pm\frac{1}{2}}w\big)^{-1}\cdot\epsilon(\alpha_i,\alpha_j)^{\pm1},
& \ (\alpha_i|\alpha_j)=-1;\vspace{1mm}\\
\big(z-(r^{-1}s)^{\pm\frac{1}{2}}w\big)\cdot\epsilon(\alpha_i,\alpha_j)^{\pm1},
& \ (\alpha_i|\alpha_j)=1;\vspace{1mm}\\
\big(z-w\big)\big(z-(r^{-1}s)^{\pm1}w\big)\cdot\epsilon(\alpha_i,\alpha_j)^{\pm1},
& \ (\alpha_i|\alpha_j)=2.\vspace{1mm}\\
\end{array}\right.
\label{OPE:1}  \\
&&Y_i^{\pm}(z)Y_j^{\mp}(w)=:Y_i^{\pm}(z)Y_j^{\mp}(w): \nonumber \\
&&\hskip2.5cm \cdot \left\{\begin{array}{lr}1, &\ (\alpha_i|\alpha_j)=0;\vspace{1mm}\\
\big(z-(rs)^{\frac{1}{4}}(rs)^{\mp\frac{1}{2}}w\big)_{r,s}^{\frac{1}{2}}\cdot\epsilon(\alpha_i,\alpha_j)^{\mp1},
& \ (\alpha_i|\alpha_j)=-\frac{1}{2};\vspace{1mm}\\
\big(z-(rs)^{\mp\frac{1}{2}}w\big)\cdot\epsilon(\alpha_i,\alpha_j)^{\mp1},
& \ (\alpha_i|\alpha_j)=-1;\vspace{1mm}\\
\big(z-(rs)^{\mp\frac{1}{2}}w\big)^{-1}\cdot\epsilon(\alpha_i,\alpha_j)^{\mp1},
& \ (\alpha_i|\alpha_j)=1;\vspace{1mm}\\
\big(z-w\big)^{-1}\big(z-(rs)^{\mp1}w\big)^{-1}\cdot\epsilon(\alpha_i,\alpha_j)^{\mp1},
& \ (\alpha_i|\alpha_j)=2.\vspace{1mm}\\
\end{array}\right.
\label{OPE:2}   \\
&&U_i^{\pm}(z)U_j^{\pm}(w)=:U_i^{\pm}(z)U_j^{\pm}(w): \nonumber\\
&&\hskip3.6cm \cdot \left\{\begin{array}{lr}1, &\ (\alpha_i|\alpha_j)=0;\vspace{1mm}\\
\big(z-(rs)^{\frac{3}{4}}w\big)_{r,s}^{-\frac{1}{2}},
& \ (\alpha_i|\alpha_j)=-\frac{1}{2};\vspace{1mm}\\
\big(z-w\big)^{-1},
& \ (\alpha_i|\alpha_j)=-1;\vspace{1mm}\\
\big(z-w\big),
& \ (\alpha_i|\alpha_j)=1;\vspace{1mm}\\
\big(z-(rs^{-1})^{\frac{1}{2}}w\big)\big(z-(r^{-1}s)^{\frac{1}{2}}w\big),
& \ (\alpha_i|\alpha_j)=2.\vspace{1mm}\\
\end{array}\right.  \label{OPE:3}\\
&&U_i^{\pm}(z)U_j^{\mp}(w)=:U_i^{\pm}(z)U_j^{\mp}(w): \nonumber\\
&&\hskip2.8cm \cdot \left\{\begin{array}{lr}1, &\ (\alpha_i|\alpha_j)=0;\vspace{1mm}\\
\big(z-(rs)^{\frac{1}{4}}w\big)_{r,s}^{\frac{1}{2}},
& \ (\alpha_i|\alpha_j)=-\frac{1}{2};\vspace{1mm}\\
\big(z-w\big),
& \ (\alpha_i|\alpha_j)=-1;\vspace{1mm}\\
\big(z-w\big)^{-1},
& \ (\alpha_i|\alpha_j)=1;\vspace{1mm}\\
\big(z-(rs^{-1})^{\frac{1}{2}}w\big)^{-1}\big(z-(r^{-1}s)^{\frac{1}{2}}w\big)^{-1},
& \ (\alpha_i|\alpha_j)=2.\vspace{1mm}\\
\end{array}\right.  \label{OPE:4}
\end{eqnarray}}
\end{lemm}
\begin{proof}\, (\ref{OPE:1}) follows from the following
\begin{eqnarray*}
\begin{split}
&LHS=Y_i^{\pm}(z)Y_j^{\pm}(w)\\
&\quad=\exp (\pm  \sum^{\infty}_{k=1}
                 \frac{a_i(-k)}{[k]} s^{\pm\frac{k}{2}} z^k)\times\exp ( \mp \sum^{\infty}_{k=1}
           \frac{a_i(k)   }{[k]} r^{\mp \frac{k}{2}} z^{-k})
         e^{\pm\a_i} z^{\pm a_i(0)}\vep_i \\
&\qquad\exp (\pm  \sum^{\infty}_{k=1}
                 \frac{a_j(-k)}{[k]} s^{\pm\frac{k}{2}} w^k)\times\exp ( \mp \sum^{\infty}_{k=1}
           \frac{a_j(k)   }{[k]} r^{\mp \frac{k}{2}} w^{-k})
         e^{\pm\a_j} z^{\pm a_j(0)}\vep_j\\
&\quad=\epsilon(\alpha_i,\alpha_j)^{\pm1}:Y_i^{\pm}(z)Y_j^{\pm}(w):\exp (-\sum^{\infty}_{k=1}
           \frac{(r^{-1}s)^{\pm \frac{k}{2}}}{[k]^2}[\,a_i(k), a_j(-k)\,] (\frac{w}{z})^{k})z^{(\alpha_i|\alpha_j)}\\
           &\quad=RHS.
\end{split}
\end{eqnarray*}
(\ref{OPE:2}) follows from the following
\begin{eqnarray*}
\begin{split}
&LHS=Y_i^{\pm}(z)Y_j^{\mp}(w)\\
&\quad=\exp (\pm  \sum^{\infty}_{k=1}
                 \frac{a_i(-k)}{[k]} s^{\pm\frac{k}{2}} z^k)\times\exp ( \mp \sum^{\infty}_{k=1}
           \frac{a_i(k)   }{[k]} r^{\mp \frac{k}{2}} z^{-k})
         e^{\pm\a_i} z^{\pm a_i(0)}\vep_i \\
&\qquad\exp (\mp  \sum^{\infty}_{k=1}
                 \frac{a_j(-k)}{[k]} s^{\mp\frac{k}{2}} w^k)\times\exp ( \pm \sum^{\infty}_{k=1}
           \frac{a_j(k)   }{[k]} r^{\pm \frac{k}{2}} w^{-k})
         e^{\mp\a_j} z^{\mp a_j(0)}\vep_j\\
&\quad=\epsilon(\alpha_i,\alpha_j)^{\mp1}:Y_i^{\pm}(z)Y_j^{\mp}(w):\exp (\sum^{\infty}_{k=1}
           \frac{(rs)^{\mp \frac{k}{2}}}{[k]^2}[\,a_i(k), a_j(-k)\,] (\frac{w}{z})^{k})z^{-(\alpha_i|\alpha_j)}\\
           &\quad=RHS.
\end{split}
\end{eqnarray*}

(\ref{OPE:3}) follows from the following
\begin{eqnarray*}
\begin{split}
&LHS=U_i^{\pm}(z)U_j^{\pm}(w)\\
&\quad=\exp (\pm  \sum^{\infty}_{k=1}
                 \frac{b_i(-k)}{[k]} z^k)\times\exp ( \mp \sum^{\infty}_{k=1}
           \frac{b_i(k)   }{[k]} z^{-k})
         e^{\pm\a_i} z^{\pm a_i(0)}(rs)^{-\frac{1}{8}} \\
&\qquad\exp (\pm  \sum^{\infty}_{k=1}
                 \frac{b_j(-k)}{[k]} s^{\pm\frac{k}{2}} w^k)\times\exp ( \mp \sum^{\infty}_{k=1}
           \frac{b_j(k)   }{[k]} r^{\mp \frac{k}{2}} w^{-k})
         e^{\pm\a_j} z^{\pm a_j(0)}(rs)^{-\frac{1}{8}} \\
&\quad=:U_i^{\pm}(z)U_j^{\pm}(w):\exp (-\sum^{\infty}_{k=1}
           \frac{1}{[k]^2}[\,b_i(k), b_j(-k)\,] (\frac{w}{z})^{k})z^{(\alpha_i|\alpha_j)}=RHS.\\
\end{split}
\end{eqnarray*}

(\ref{OPE:4}) follows from the following
\begin{eqnarray*}
\begin{split}
&LHS=U_i^{\pm}(z)U_j^{\mp}(w)\\
&\quad=\exp (\pm  \sum^{\infty}_{k=1}
                 \frac{b_i(-k)}{[k]} z^k)\times\exp ( \mp \sum^{\infty}_{k=1}
           \frac{b_i(k)   }{[k]} z^{-k})
         e^{\pm\a_i} z^{\pm a_i(0)}(rs)^{-\frac{1}{8}} \\
&\qquad\exp (\mp  \sum^{\infty}_{k=1}
                 \frac{b_j(-k)}{[k]} s^{\mp\frac{k}{2}} w^k)\times\exp ( \pm \sum^{\infty}_{k=1}
           \frac{b_j(k)   }{[k]} r^{\pm \frac{k}{2}} w^{-k})
         e^{\mp\a_j} z^{\mp a_j(0)}(rs)^{-\frac{1}{8}} \\
&\quad=:U_i^{\pm}(z)U_j^{\mp}(w):\exp (\sum^{\infty}_{k=1}
           \frac{1}{[k]^2}[\,b_i(k), b_j(-k)\,] (\frac{w}{z})^{k})z^{-(\alpha_i|\alpha_j)}=RHS.
\end{split}
\end{eqnarray*}
\end{proof}

Now we proceed to check relation $(D7)$, which holds from the following Proposition.
\begin{prop}
\begin{eqnarray*}
&&(z-(\lg i,j\rg\lg
j,i\rg)^{\pm\frac1{2}}w)\,X_i^{\pm}(z)X_j^{\pm}(w)\\
&&\qquad =(\lg
j,i\rg^{\pm1}z-(\lg j,i\rg\lg
i,j\rg^{-1})^{\pm\frac1{2}}w)\,X_j^{\pm}(w)\,X_i^{\pm}(z).
\end{eqnarray*}
\end{prop}
\begin{proof} \ The proof will be carried out in the following two cases, the other cases can be checked similarly.

 $(\text{i})$ For the case of $(\alpha_i|\alpha_j)=-1$, let us first consider
 \begin{eqnarray*}
 \begin{split}
 &X^-_{n-1}(z)X^-_{n}(w)=:X^-_{n-1}(z)X^-_{n}(w): \big(z-(r^{-1}s)^{-\frac{1}{2}}w\big)^{-1}.
\end{split}
\end{eqnarray*}
On the other hand, one gets
\begin{eqnarray*}
 \begin{split}
 & X^-_{n}(w)X^-_{n-1}(z)=-(rs)^{\frac{1}{2}}:X^-_{n}(w)X^-_{n-1}(z): \big(w-(r^{-1}s)^{-\frac{1}{2}}z\big)^{-1}.
\end{split}
\end{eqnarray*}
It holds that
$$(z-(rs^{-1})^{\frac{1}{2}}w)X^-_{n-1}(z)X^-_{n}(w)
=(s^{-1}z-(rs)^{-\frac{1}{2}}w)X^-_{n}(w)X^-_{n-1}(z).$$

$(\text{ii})$ For the case of $(\alpha_i|\alpha_j)=-\frac{1}{2}$, that is, $(\alpha_i|\alpha_{i+1})=-\frac{1}{2}$.
 For $i=1, \cdots, n-2$, by (4.14)--(4.17), one has
  \begin{eqnarray*}
&&X^+_i(z)X^+_{i+1}(w)\\
&&= \varepsilon(\alpha_{i},\,\alpha_{i+1}):Y_i^+(z)Y_{i+1}^+(w):\big(z-(rs)^{\frac{3}{4}}(r^{-1}s)^{\frac{1}{2}}w\big)_{r,s}^{-\frac{1}{2}}\\
 && \left(:U_i^+(zs^{\frac{1}{2}})U^+_{i+1}(wr^{-\frac{1}{2}}):\big(z-(rs)^{\frac{3}{4}}w\big)_{r,s}^{-\frac{1}{2}}
 \cdot s^{\frac{1}{2}(\alpha_i|\alpha_{i+1})} \right.\\
&& +:U_i^+(zs^{\frac{1}{2}})U_{i+1}^-(wr^{\frac{1}{2}}):(-1)^{2a_{i+1}(0)} \big(z-(rs)^{\frac{1}{4}}(rs^{-1})^{\frac{1}{2}}w\big)_{r,s}^{\frac{1}{2}}\cdot s^{-\frac{1}{2}(\alpha_i|\alpha_{i+1})}\\
 &&  +:U_i^-(zr^{\frac{1}{2}})U_{i+1}^+(ws^{\frac{1}{2}}):(-1)^{2a_{i}(0)+1}
 \big(z-(rs)^{\frac{1}{4}}(r^{-1}s)^{\frac{1}{2}}w\big)_{r,s}^{\frac{1}{2}}\cdot r^{-\frac{1}{2}(\alpha_i|\alpha_{i+1})}\\
&&  \left.+
:U_i^-(zr^{\frac{1}{2}})U_{i+1}^-(wr^{\frac{1}{2}}):(-1)^{2a_i(0)+2a_{i+1}(0)+1}
\big(z-(rs)^{\frac{3}{4}}w\big)_{r,s}^{-\frac{1}{2}}\cdot r^{\frac{1}{2}(\alpha_i|\alpha_{i+1})} \right).
\end{eqnarray*}
Using (4.9)--(4.11), one easily gets
\begin{eqnarray*}
&&X^+_i(z)X^+_{i+1}(w)\\
&=& :Y_i^+Y_{i+1}^+:
\big(:U_i^+(s^{\frac{1}{2}}z)U^+_{i+1}(s^{\frac{1}{2}}w):(z-(r^{-1}s)^{\frac{1}{4}}w)^{-1}s^{-\frac{1}{4}}\\
&&+:U_{i}^+(s^{\frac{1}{2}}z)U_{i+1}^-(r^{\frac{1}{2}}w):(-1)^{2a_{i+1}(0)}
\frac{z-(rs^{-1})^{\frac{1}{4}}w}{z-(r^{-1}s)^{\frac{1}{4}}w}s^{\frac{1}{4}}\\
&&+:U_i^-(r^{\frac{1}{2}}z)U^+_{i+1}(s^{\frac{1}{2}}w):(-1)^{2a_{i}(0)+1}r^{\frac{1}{4}}\\
&&+ :U_i^-(r^{\frac{1}{2}}z)U_{i+1}^-(r^{\frac{1}{2}}w):(-1)^{2a_i(0)+2a_{i+1}(0)+1}
(z-(r^{-1}s)^{\frac{1}{4}}w)^{-1}r^{-\frac{1}{4}}\big).
\end{eqnarray*}
On the other hand, it holds similarly
\begin{eqnarray*}
&&X^+_{i+1}(w)X^+_i(z)\\
&=& -(rs)^{-\frac{1}{4}}:Y_{i+1}^+Y_i^+:
\big(:U^+_{i+1}(s^{\frac{1}{2}}w)U^+_i(s^{\frac{1}{2}}z):(w-(r^{-1}s)^{\frac{1}{4}}z)^{-1}\cdot s^{-\frac{1}{4}}\\
&&+:U_{i+1}^-(r^{\frac{1}{2}}w)U_{i}^+(s^{\frac{1}{2}}z):(-1)^{2a_{i+1}(0)+1}\cdot r^{\frac{1}{4}}\\
&&+:U^+_{i+1}(s^{\frac{1}{2}}w)U_i^-(r^{\frac{1}{2}}z):
(-1)^{2a_{i}(0)}\frac{w-(rs^{-1})^{\frac{1}{4}}z}{w-(r^{-1}s)^{\frac{1}{4}}z}\cdot s^{\frac{1}{4}}\\
&&+ :U_{i+1}^-(r^{\frac{1}{2}}w)U_i^-(r^{\frac{1}{2}}z):(-1)^{2a_i(0)+2a_{i+1}(0)+1}
(w-(r^{-1}s)^{\frac{1}{4}}z)^{-1}\cdot r^{-\frac{1}{4}}\big).
\end{eqnarray*}

As a consequence of the above relations, we get our required result
\begin{eqnarray*}
&&(z-(r^{-1}s)^{\frac{1}{4}}w)X_i^+(z)X_{i+1}^+(w)
=(s^{\frac{1}{2}}z-(rs)^{\frac{1}{4}}w)X_{i+1}^+(w)X_{i}^+(z).
\end{eqnarray*}

The ``-"-part of relation $({\rm D}7)$ can be verified similarly, which is left to the readers.
\end{proof}

Now we focus on checking relation $({\rm D}8)$.
\begin{lemm}
$$[\,X_i^+(z),
X_j^-(w)\,]=\frac{\delta_{ij}}{(r_i-s_i)zw}\Big(\delta(zw^{-1}s)\psi_i(wr^{\frac{1}2})
-\delta(zw^{-1}r)\phi_i(ws^{-\frac1{2}})\Big).$$
\end{lemm}
\begin{proof} It suffices to consider the four
cases: $(\a_i|\a_j)=-1/2$, $(\a_i|\a_j)=1$,  $(\a_i|\a_j)=-1$, and
$(\a_i|\a_j)=2$. Here we only give the proof for the first two cases,
since the other cases are either immediate or similar to our
previous considerations.

(i) For $(\a_i|\a_j)=-1/2$, that is, $j=i+1$ or $i-1$. Firstly, consider $i=1,\cdots, n-1$,
using (4.14)--(4.17), together with (4.11)--(4.13), one immediately gets
\begin{eqnarray*}
&&X^+_i(z)X_{i+1}^-(w)=:Y_i^+(z)Y_{i+1}^-(w):\\
&&\left(:U_i^+(s^{\frac{1}{2}}z)U_{i+1}^+(r^{-\frac{1}{2}}w):s^{-\frac{1}{4}}\right. \\
&&+:U_i^+(s^{\frac{1}{2}}z)U_{i+1}^-(s^{-\frac{1}{2}}w):(-1)^{2a_{i+1}(0)}
(s^{\frac{1}{4}}z-r^{-\frac{1}{4}}s^{-\frac{1}{2}}w)\\
&&+:U_i^-(r^{\frac{1}{2}}z)U_i^+(r^{-\frac{1}{2}}w):(-1)^{2a_i(0)+1}
(r^{\frac{1}{4}}z-r^{-\frac{1}{2}}s^{-\frac{1}{4}}w) \\
&&+\left.:U_i^-(r^{\frac{1}{2}}z)U_{i+1}^-(s^{-\frac{1}{2}}w):(-1)^{2a_i(0)+2a_{i+1}(0)+1}r^{-\frac{1}{4}}\right).
\end{eqnarray*}
Using the same method, it is easy to see that
\begin{eqnarray*}
&&X_{i+1}^-(w)X^+_i(z)=(rs)^{-\frac{1}{4}}
:Y_{i+1}^-(w)Y_i^+(z):\epsilon_0(\alpha_{i+1},\,\alpha_{i})\\
&&\left(:U_{i+1}^+(r^{-\frac{1}{2}}w)U_i^+(s^{\frac{1}{2}}z):r^{\frac{1}{4}}\right. \\
&&+:U_{i+1}^-(s^{-\frac{1}{2}}w)U_i^+(s^{\frac{1}{2}}z):(-1)^{2a_{i+1}(0)+1}
(s^{-\frac{1}{4}}z-r^{\frac{1}{4}}s^{\frac{1}{2}}w) \\
&&+:U_{i+1}^+(r^{-\frac{1}{2}}w)U_i^-(r^{\frac{1}{2}}z):(-1)^{2a_i(0)}
(r^{-\frac{1}{4}}z-r^{\frac{1}{2}}s^{\frac{1}{4}}w) \\
&&+\left.:U_{i+1}^-(s^{-\frac{1}{2}}w)U_i^-(r^{\frac{1}{2}}z):(-1)^{4a_i(0)}s^{\frac{1}{4}}\right).
\end{eqnarray*}
Thus it is easy to get  $[X^+_{i}(z), X^-_{i+1}(w)]=0.$
\smallskip

(ii)\ For the case of $1\leqslant i=j\leqslant n-1$,
using (\ref{OPE:1})-(\ref{OPE:4}), it follows from (4.11)--(4.13) that
\begin{eqnarray*}
&&X^+_i(z)X_i^-(w)=(rs)^{\frac{1}{4}}:Y_i^+(z)Y_i^-(w):\\
&&\left(:U_i^+(s^{\frac{1}{2}}z)U_i^+(r^{-\frac{1}{2}}w):s^{\frac{1}{2}}\right. \\
&&+:U_i^+(s^{\frac{1}{2}}z)U_i^-(s^{-\frac{1}{2}}w):(-1)^{2a_i(0)}\frac
1{(s^{\frac{1}{2}}z-s^{-\frac{1}{2}}w)
(z-(rs)^{-\frac{1}{2}}w)} \\
&&+:U_i^-(r^{\frac{1}{2}}z)U_i^+(r^{-\frac{1}{2}}w):(-1)^{2a_i(0)}\frac
1{(r^{\frac{1}{2}}z-r^{-\frac{1}{2}}w)
(z-(rs)^{-\frac{1}{2}}w)} \\
&&+\left.:U_i^-(r^{\frac{1}{2}}z)U_i^-(s^{-\frac{1}{2}}w):(-1)^{4a_i(0)}r^{\frac{1}{2}}\right).
\end{eqnarray*}
At the same time, we actually have
\begin{eqnarray*}
&&X_i^-(w)X^+_i(z)=(rs)^{-\frac{1}{4}}:Y_i^-(w)Y_i^+(z):\\
&&\left(:U_i^+(r^{-\frac{1}{2}}w)U_i^+(s^{\frac{1}{2}}z):r^{-\frac{1}{2}}\right. \\
&&+:U_i^-(s^{-\frac{1}{2}}w)U_i^+(s^{\frac{1}{2}}z):(-1)^{2a_i(0)}\frac
1{(s^{-\frac{1}{2}}w-s^{\frac{1}{2}}z)
(w-(rs)^{\frac{1}{2}}z)} \\
&&+:U_i^+(r^{-\frac{1}{2}}w)U_i^-(r^{\frac{1}{2}}z):(-1)^{2a_i(0)}\frac
1{(r^{-\frac{1}{2}}w-r^{\frac{1}{2}}z)
(w-(rs)^{\frac{1}{2}}z)} \\
&&+\left.:U_i^-(s^{-\frac{1}{2}}w)U_i^-(r^{\frac{1}{2}}z):(-1)^{4a_i(0)}s^{-\frac{1}{2}}\right).
\end{eqnarray*}

Therefore, it is easy to see that
\begin{eqnarray*}
&&[X^+_{i}(z), X^-_{i}(w)]\\
&&=(rs)^{-\frac{1}{4}}:U_i^+(s^{\frac{1}{2}}z)U_i^-(s^{-\frac{1}{2}}w)Y_i^+(z)Y_i^-(w):
\left(\frac1{(s^{\frac{1}{2}}z-s^{-\frac{1}{2}}w)
(z-(rs)^{-\frac{1}{2}}w)}\right.\\
&&\qquad\quad \left. -\frac 1{(s^{-\frac{1}{2}}w-s^{\frac{1}{2}}z)
((rs)^{-\frac{1}{2}}w-z)}\right)\\
&&+(rs)^{-\frac{1}{4}}:U_i^-(r^{\frac{1}{2}}z)U_i^+(r^{-\frac{1}{2}}w)Y_i^+(z)Y_i^-(w):
\left(\frac1{(r^{\frac{1}{2}}z-r^{-\frac{1}{2}}w)
(z-(rs)^{-\frac{1}{2}}w)}\right.\\
&&\qquad\quad \left.-\frac 1{(r^{-\frac{1}{2}}w-r^{\frac{1}{2}}z)
((rs)^{-\frac{1}{2}}w-z)}\right)\\
&&=:U_i^+(s^{\frac{1}{2}}z)U_i^-(s^{-\frac{1}{2}}w)Y_i^+(z)Y_i^-(w):
\frac{(rs)^{\frac{1}{4}}}{(r^{\frac{1}{2}}-s^{\frac{1}{2}})zw}
\left(\delta(s^{-1}\frac wz)-\delta(\frac{(rs)^{-\frac{1}{2}}w}{z} )
\right)\\
&&+:U_i^-(r^{\frac{1}{2}}z)U_i^+(r^{-\frac{1}{2}}w)Y_i^+(z)Y_i^-(w):
\frac{(rs)^{\frac{1}{4}}}{(r^{\frac{1}{2}}-s^{\frac{1}{2}})zw}
\left(\delta(\frac{(rs)^{-\frac{1}{2}}w}{z}-\delta(r^{-1}\frac wz)
 )\right)\\
&&=\frac1{(r_i-s_i)zw}
                 \left(
                 \psi_i(wr^{\frac{1}{2}})
                 \delta(\frac{ws^{-1}}{z})
                 -\varphi_i(ws^{-\frac{1}{2}})
                  \delta(\frac{wr^{-1}}z)
                 \right),
\end{eqnarray*}
where we have used the following results
\begin{eqnarray*}
&&:U_i^+(s^{-\frac{1}{2}}w)U_i^-(s^{-\frac{1}{2}}w)Y_i^+(s^{-1}w)Y_i^-(w):
=(rs)^{-\frac{1}{4}}\psi_i(wr^{\frac{1}{2}});\\
&&:U_i^-(r^{-\frac{1}{2}}w)U_i^+(r^{-\frac{1}{2}}w)Y_i^+(r^{-1}w)Y_i^-(w):=(rs)^{-\frac{1}{4}}\varphi_i(ws^{-\frac{1}{2}}).
\end{eqnarray*}

This completes the proof.
\end{proof}

\noindent{\bf Proof of Serre relations $({\rm D}9)$}. (i) \ For $i=1, \cdots, n-2$,
let us write the operators $X^{\pm}_i(z)$ as a sum of two terms:
\begin{eqnarray*}
X^{+}_i(z)&=&\sum_{\ep=\pm}Y_i^{+}(z)
U_i^{\ep}(zr^{\frac{(1-\ep)}{4}}s^{\frac{(\ep+1)}{4}})
(-1)^{(1-\ep)a_i(0)}\\
&=& X^{+}_{i+}(z)+X^{+}_{i-}(z),\\
X^{-}_i(z)&=&\sum_{\ep=\pm}Y_i^{-}(z)
U_i^{\ep}(zr^{-\frac{(\ep+1)}{4}}s^{\frac{(\ep-1)}{4}})
(-1)^{(1-\ep)a_i(0)}\\
&=& X^{-}_{i+}(z)+X^{-}_{i-}(z).
\end{eqnarray*}

From the relations (\ref{OPE:1})-(\ref{OPE:4}), it follows that
\begin{eqnarray*}
&&X_{i\ep_1}^+(z_1)X_{i\ep_2}^+(z_2)X^+_{i+1,\ep}(w)=
:X_{i\ep_1}^+(z_1)X_{i\ep_2}^+(z_2)X^+_{i+1,\ep}(w):(-1)^{(\ep_1-\ep_2)/2}\\
&&\qquad \times (z_1-(r^{-1}s)^{\frac{1}{2}}z_2)
(r^{\frac{1-\ep_1}{4}}s^{\frac{\ep_1+1}{4}}z_1-
r^{\frac{1-\ep_2}{4}}s^{\frac{\ep_2+1}{4}}z_2)^{\ep_1\ep_2}
r^{\frac{(2-\ep_1-\ep_2)\ep}{8}}s^{-\frac{(2+\ep_1+\ep_2)\ep}{8}}  \\
&&\qquad \times
\frac{(z_1-(r^{-1}s)^{\frac{\ep}{4}}w)^{\frac{|\ep-\ep_1|}2}}{z_1-(r^{-1}s)^{\frac{1}{4}}w}
\frac{(z_2-(r^{-1}s)^{\frac{\ep}{4}}w)^{\frac{|\ep-\ep_2|}2}}{z_2-(r^{-1}s)^{\frac{1}{4}}w},
\end{eqnarray*}
where we include the sign factor $(-1)^{(1-\ep)\a_i(0)}$ and
$\vep_i$ in the normal ordered product. Similar normal product
computation gives that
\begin{eqnarray*}
&&X_{i\ep_1}^+(z_1)X^+_{i+1,\ep}(w)X_{i\ep_2}^+(z_2)=
:X_{i\ep_1}^+(z_1)X^+_{i+1,\ep}(w)X_{i\ep_2}^+(z_2):(-1)^{(\ep_1-\ep_2)/2}\\
&&\qquad \times (z_1-(r^{-1}s)^{\frac{1}{2}}z_2)
(r^{\frac{1-\ep_1}{4}}s^{\frac{\ep_1+1}{4}}z_1-
r^{\frac{1-\ep_2}{4}}s^{\frac{\ep_2+1}{4}}z_2)^{\ep_1\ep_2}
r^{\frac{(2-\ep_1-\ep_2)\ep}{8}}s^{-\frac{(2+\ep_1+\ep_2)\ep}{8}}\\
&&\qquad \times
\frac{(z_1-(r^{-1}s)^{\frac{\ep}{4}}w)^{\frac{|\ep-\ep_1|}2}}{z_1-(r^{-1}s)^{\frac{1}{4}}w}
\frac{((r^{-1}s)^{\frac{\ep}{4}}w-z_2)^{\frac{|\ep-\ep_2|}2}}{w-(r^{-1}s)^{\frac{1}{4}}z_2}.
\end{eqnarray*}
\begin{eqnarray*}
&&X^+_{i+1,\ep}(w)X_{i\ep_1}^+(z_1)X_{i\ep_2}^+(z_2)=
:X^+_{i+1,\ep}(w)X_{i\ep_1}^+(z_1)X_{i\ep_2}^+(z_2):(-1)^{(\ep_1-\ep_2)/2}\\
&&\qquad \times (z_1-(r^{-1}s)^{\frac{1}{2}}z_2)
(r^{\frac{1-\ep_1}{4}}s^{\frac{\ep_1+1}{4}}z_1-
r^{\frac{1-\ep_2}{4}}s^{\frac{\ep_2+1}{4}}z_2)^{\ep_1\ep_2}
r^{\frac{(2-\ep_1-\ep_2)\ep}{8}}s^{-\frac{(2+\ep_1+\ep_2)\ep}{8}} \\
&&\qquad \times
\frac{((r^{-1}s)^{\frac{\ep}{4}}w-z_1)^{\frac{|\ep-\ep_1|}2}}{w-(r^{-1}s)^{\frac{1}{4}}z_1}
\frac{((r^{-1}s)^{\frac{\ep}{4}}w-z_2)^{\frac{|\ep-\ep_2|}2}}{w-(r^{-1}s)^{\frac{1}{4}}z_2}.
\end{eqnarray*}
One may now use above results to get that
\begin{eqnarray*}
&&X_{i\ep_1}^+(z_1)X_{i\ep_2}^+(z_2)X^+_{i+1,\ep}(w)-(r^{\frac{1}{2}}+s^{\frac{1}{2}})
X_{i\ep_1}^+(z_1)X^+_{i+1,\ep}(w)X_{i\ep_2}^+(z_2)\\
&&\quad +(rs)^{\frac{1}{2}}X^+_{i+1,\ep}(w)X_{i\ep_1}^+(z_1)X_{i\ep_2}^+(z_2) \\
&=&:X_{i\ep_1}^+(z_1)X_{i\ep_2}^+(z_2)X^+_{i+1,\ep}(w):
(z_1-(r^{-1}s)^{\frac{1}{2}}z_2)
(r^{\frac{1-\ep_1}{4}}s^{\frac{\ep_1+1}{4}}z_1-
r^{\frac{1-\ep_2}{4}}s^{\frac{\ep_2+1}{4}}z_2)^{\ep_1\ep_2}\\
&& \quad \cdot (-1)^{(\ep_1-\ep_2)/2}\times r^{\frac{(2-\ep_1-\ep_2)\ep}{8}}s^{-\frac{(2+\ep_1+\ep_2)\ep}{8}}
\Big( \frac{(z_1-(r^{-1}s)^{\frac{\ep}{4}}w)^{\frac{|\ep-\ep_1|}2}}{z_1-(r^{-1}s)^{\frac{1}{4}}w}\times
 \end{eqnarray*}
\begin{eqnarray} \label{g:1}
&&\quad\frac{(z_2-(r^{-1}s)^{\frac{\ep}{4}}w)^{\frac{|\ep-\ep_2|}2}}{z_2-(r^{-1}s)^{\frac{1}{4}}w}
 +[2]_{\frac{1}{2}}(rs)^{-\frac{1}{4}}\times
\frac{(z_1-(r^{-1}s)^{\frac{\ep}{4}}w)^{\frac{|\ep-\ep_1|}2}}{z_1-(r^{-1}s)^{\frac{1}{4}}w}\times \\
&&\quad\frac{((r^{-1}s)^{\frac{\ep}{4}}w-z_2)^{\frac{|\ep-\ep_2|}2}}{w-(r^{-1}s)^{\frac{1}{4}}z_2}+
\frac{((r^{-1}s)^{\frac{\ep}{4}}w-z_1)^{\frac{|\ep-\ep_1|}2}}{w-(r^{-1}s)^{\frac{1}{4}}z_1}
\frac{((r^{-1}s)^{\frac{\ep}{4}}w-z_2)^{\frac{|\ep-\ep_2|}2}}{w-(r^{-1}s)^{\frac{1}{4}}z_2}\Big),\nonumber
  \end{eqnarray}
where each term in the parentheses corresponds to the above relations for the
three normal products, and we have also used the relation:
$$:X_{i\ep_1}^+(z_1)X_{i\ep_2}^+(z_2)X^+_{i+1,\ep}(w):
=-(rs)^{\frac{1}{2}}:X_{i\ep_1}^+(z_1)X^+_{i+1,\ep}(w)X_{i\ep_2}^+(z_2):.$$

We check the claim for two cases $\ep_1\neq \ep_2$ and $\ep_1=\ep_2$ :

(a) \ For $\ep_1\neq \ep_2$, one has
\begin{eqnarray}\label{E:S2}
&&X_{i\ep_1}^+(z_1)X_{i\ep_2}^+(z_2)X^+_{i+1,\ep}(w)-(r^{\frac{1}{2}}+s^{\frac{1}{2}})
X_{i\ep_1}^+(z_1)X^+_{i+1,\ep}(w)X_{i\ep_2}^+(z_2)\nonumber\\
&&\quad
+(rs)^{\frac{1}{2}}X^+_{i+1,\ep}(w)X_{i\ep_1}^+(z_1)X_{i\ep_2}^+(z_2)
+(z_1\leftrightarrow z_2, \ep_1\leftrightarrow \ep_2)=0.
\end{eqnarray}

The claim of this case is verified by checking four cases for $\ep, \ep_i$, which
are all similar and relied upon the following important identity
\cite{JKM1}:
\begin{eqnarray} \label{identity1}
&&(z_1-tw)(z_2-tw)
+(t+t^{-1})(z_1-tw)
(w-tz_2)%\nonumber\\
%&&\qquad\qquad
+(w-tz_1)
(w-tz_2)\\
&&=(t^{-1}-t)w(z_1-t^{2}z_2),\nonumber
\end{eqnarray}
for any  $t\in {\bf C}$.

Take $\ep=1, \ep_1=-\ep_2=1$ for example, the parentheses in (\ref{g:1})
is simplified to the following expression times
$\prod_i(z_i-(r^{-1}s)^{\frac{1}{4}}w)^{-1} \cdot
(w-(r^{-1}s)^{\frac{1}{4}}z_2)^{-1}$.
\begin{eqnarray*}
&&\left((z_1-(r^{-1}s)^{\frac{1}{4}}w)(z_2-(r^{-1}s)^{\frac{1}{4}}w)+(rs)^{-\frac{1}{4}}[2]_i
(z_1-(r^{-1}s)^{\frac{1}{4}}w)(w-(r^{-1}s)^{\frac{1}{4}}z_2)
\right.
\\
&&\quad \left. +(w-(r^{-1}s)^{\frac{1}{4}}z_1)(w-(r^{-1}s)^{\frac{1}{4}}z_2)\right)\\
&&=((rs^{-1})^{\frac{1}{4}}-(r^{-1}s)^{\frac{1}{4}})w(z_1-(r^{-1}s)^{\frac{1}{2}}z_2).
\end{eqnarray*}

Under the symmetry $(z_1, \ep_1)\leftrightarrow (z_2, \ep_2)$ it
follows that the claim holds due to
$$w(z_1-(r^{-1}s)^{\frac{1}{2}}z_2)+w((r^{-1}s)^{\frac{1}{2}}z_2-z_1)=0.$$

(b) \ We now turn to the other four cases with $\ep_1=\ep_2$. For the case $\ep=\ep_1=\ep_2=-1$, using
the identity (\ref{identity1}) again to simplify the parentheses in
(\ref{g:1}), the contraction function in the Serre relation becomes
\begin{eqnarray*}
&&s(z_1-(r^{-1}s)^{\frac{1}{2}}z_2)(z_1-z_2)
\left(\frac{(z_1-(rs^{-1})^{\frac{1}{4}}w)(z_2-(rs^{-1})^{\frac{1}{4}}w)}
{(z_1-(r^{-1}s)^{\frac{1}{4}}w)(z_2-(r^{-1}s)^{\frac{1}{4}}w)}
\right.\\
&&\qquad\quad \left.-[2]_{\frac{1}{2}}(rs)^{-\frac{1}{4}}
\frac{z_1-(rs^{-1})^{\frac{1}{4}}w}{z_1-(r^{-1}s)^{\frac{1}{4}}w}+(rs^{-1})^{\frac{1}{4}}\right)\\
&=& \frac{s((r^{-1}s)^{\frac{1}{4}}-(rs^{-1})^{\frac{1}{4}})w
(z_1-z_2)(z_1-(r^{-1}s)^{\frac{1}{2}}z_2)(z_1-(rs^{-1})^{\frac{1}{2}}z_2)}
{(z_1-(r^{-1}s)^{\frac{1}{4}}w)(z_2-(r^{-1}s)^{\frac{1}{4}}w)},
\end{eqnarray*}
which is anti-symmetric under $(z_1\leftrightarrow z_2)$, hence the
sub-Serre relation is proved in this case. That is,
\begin{eqnarray*}
&&X_{i\ep_1}^+(z_1)X_{i\ep_1}^+(z_2)X^+_{i+1,\ep}(w)-(r^{\frac{1}{2}}+s^{\frac{1}{2}})
X_{i\ep_1}^+(z_1)X^+_{i+1,\ep}(w)X_{i\ep_1}^+(z_2)\\
&&\quad
+(rs)^{\frac{1}{2}}X^+_{i+1,\ep}(w)X_{i\ep_1}^+(z_1)X_{i\ep_1}^+(z_2)
+(z_1\leftrightarrow z_2)=0.
\end{eqnarray*}
Combining this sub-Serre relation with (\ref{E:S2}), we prove the
Serre relation for $a_{i, i+1}=a_{i+1. i}=-1$.

We remark that the case $a_{n-1, n}=-1$ is easily proved by using
the identity (\ref{identity1}) with $t=(r^{-1}s)^{\frac{1}{2}}$.

Finally, let us show the fourth order Serre relation with $a_{n,
n-1}=-2$.

\begin{eqnarray}
&&Sym_{z_1,
z_2,z_3}\Big(X_{n-1}^{+}(z_1)
X_{n-1}^{+}(z_2)X_{n-1}^{+}(z_3)X_n^{+}(w)
\nonumber\\
&\quad&-
(r{+}(rs)^{\frac{1}{2}}+s)\,X_{n-1}^{+}(z_1)X_{n-1}^{+}(z_2)X_n^{+}(w)X_{n-1}^{+}(z_3)\\
&\quad&+(rs)^{\frac{1}{2}}(r{+}(rs)^{\frac{1}{2}}+s)
X_{n-1}^{+}(z_1)X_n^{+}(w)X_{n-1}^{+}(z_2)X_{n-1}^{+}(z_3)\nonumber\\
&\quad&-(rs)^{\frac{3}{2}}X_n^{+}(w)X_{n-1}^{+}(z_1)X_{n-1}^{+}(z_2)X_{n-1}^{+}(z_3)\Big)=0.\nonumber
\end{eqnarray}

Firstly, using the relations (\ref{OPE:1})-(\ref{OPE:4}), it is not difficult to get
\begin{eqnarray*}
&&X_{n-1,\,\ep_1}^+(z_1)X_{n-1,\,\ep_2}^+(z_2)X_{n-1,\,\ep_3}^+(z_3)X^+_{n}(w)\\
&&=:X_{n-1,\,\ep_1}^+(z_1)X_{n-1,\,\ep_2}^+(z_2)X_{n-1,\,\ep_3}^+(z_3)X^+_{n}(w):\\
&&\quad \frac{\prod_{1\leqslant i<j\leqslant 3}(z_i-(r^{-1}s)^{\frac{1}{2}}z_j)
(r^{\frac{1-\ep_i}{4}}s^{\frac{1+\ep_i}{4}}z_i-r^{\frac{1-\ep_j}{4}}s{\frac{1+\ep_j}{4}}z_j)^{\ep_i\ep_j}}
{\prod_{i=1}^{3}(z_i-(r^{-1}s)^{\frac{1}{2}}w)
(r^{\frac{1-\ep_i}{4}}s^{\frac{1+\ep_i}{4}}z_i-r^{\frac{1-\ep_j}{4}}s{\frac{1+\ep_j}{4}}z_j)^{\ep_i}}.
\end{eqnarray*}

Similarly, it follows from the relations (\ref{OPE:1})-(\ref{OPE:4}),
\begin{eqnarray*}
&&X_{n-1,\,\ep_1}^+(z_1)X_{n-1,\,\ep_2}^+(z_2)X^+_{n}(w)X_{n-1,\,\ep_3}^+(z_3)\\
&&=:X_{n-1,\,\ep_1}^+(z_1)X_{n-1,\,\ep_2}^+(z_2)X^+_{n}(w)X_{n-1,\,\ep_3}^+(z_3):\\
&&\quad \frac{(z_3-(r^{-1}s)^{\frac{1}{2}}w)\prod_{1\leqslant i<j\leqslant 3}(z_i-(r^{-1}s)^{\frac{1}{2}}z_j)
(r^{\frac{1-\ep_i}{4}}s^{\frac{1+\ep_i}{4}}z_i-r^{\frac{1-\ep_j}{4}}s{\frac{1+\ep_j}{4}}z_j)^{\ep_i\ep_j}}
{(w-(r^{-1}s)^{\frac{1}{2}}z_3)
\prod_{i=1}^{3}(z_i-(r^{-1}s)^{\frac{1}{2}}w)
(r^{\frac{1-\ep_i}{4}}s^{\frac{1+\ep_i}{4}}z_i-r^{\frac{1-\ep_j}{4}}s{\frac{1+\ep_j}{4}}z_j)^{\ep_i}},
\end{eqnarray*}
\begin{eqnarray*}
&&X_{n-1,\,\ep_1}^+(z_1)X^+_{n}(w)X_{n-1,\,\ep_2}^+(z_2)X_{n-1,\,\ep_3}^+(z_3)\\
&&=:X_{n-1,\,\ep_1}^+(z_1)X^+_{n}(w)X_{n-1,\,\ep_2}^+(z_2)X_{n-1,\,\ep_3}^+(z_3):\\
&&\quad \frac{(w-(r^{-1}s)^{\frac{1}{2}}z_1)\prod_{1\leqslant i<j\leqslant 3}(z_i-(r^{-1}s)^{\frac{1}{2}}z_j)
(r^{\frac{1-\ep_i}{4}}s^{\frac{1+\ep_i}{4}}z_i-r^{\frac{1-\ep_j}{4}}s{\frac{1+\ep_j}{4}}z_j)^{\ep_i\ep_j}}
{(z_1-(r^{-1}s)^{\frac{1}{2}}w)\prod_{i=1}^{3}(w-(r^{-1}s)^{\frac{1}{2}}z_i)
(r^{\frac{1-\ep_i}{4}}s^{\frac{1+\ep_i}{4}}z_i-r^{\frac{1-\ep_j}{4}}s{\frac{1+\ep_j}{4}}z_j)^{\ep_i}},
\end{eqnarray*}
\begin{eqnarray*}
&&X^+_{n}(w)X_{n-1,\,\ep_1}^+(z_1)X_{n-1,\,\ep_2}^+(z_2)X_{n-1,\,\ep_3}^+(z_3)\\
&&=:X^+_{n}(w)X_{n-1,\,\ep_1}^+(z_1)X_{n-1,\,\ep_2}^+(z_2)X_{n-1,\,\ep_3}^+(z_3):\\
&&\quad \frac{\prod_{1\leqslant i<j\leqslant 3}(z_i-(r^{-1}s)^{\frac{1}{2}}z_j)
(r^{\frac{1-\ep_i}{4}}s^{\frac{1+\ep_i}{4}}z_i-r^{\frac{1-\ep_j}{4}}s{\frac{1+\ep_j}{4}}z_j)^{\ep_i\ep_j}}
{\prod_{i=1}^{3}(w-(r^{-1}s)^{\frac{1}{2}}z_i)
(r^{\frac{1-\ep_i}{4}}s^{\frac{1+\ep_i}{4}}z_i-r^{\frac{1-\ep_j}{4}}s{\frac{1+\ep_j}{4}}z_j)^{\ep_i}}.
\end{eqnarray*}

As a consequence, pulling out the common normal product, (4.19) becomes
\begin{eqnarray*}
&&Sym_{z_1,
z_2,z_3}\frac{\prod_{1\leqslant i<j\leqslant 3}(z_i-(r^{-1}s)^{\frac{1}{2}}z_j)
(r^{\frac{1-\ep_i}{4}}s^{\frac{1+\ep_i}{4}}z_i-r^{\frac{1-\ep_j}{4}}s{\frac{1+\ep_j}{4}}z_j)^{\ep_i\ep_j}}
{\prod_{i=1}^{3}(z_i-(r^{-1}s)^{\frac{1}{2}}w)(w-(r^{-1}s)^{\frac{1}{2}}z_i)
(r^{\frac{1-\ep_i}{4}}s^{\frac{1+\ep_i}{4}}z_i-r^{\frac{1-\ep_j}{4}}s{\frac{1+\ep_j}{4}}z_j)^{\ep_i}}\\
&&\Big((w-(r^{-1}s)^{\frac{1}{2}}z_1)(w-(r^{-1}s)^{\frac{1}{2}}z_2)(w-(r^{-1}s)^{\frac{1}{2}}z_3)\\
&&-(r{+}(rs)^{\frac{1}{2}}+s)\,(z_1-(r^{-1}s)^{\frac{1}{2}}w)(w-(r^{-1}s)^{\frac{1}{2}}z_2)(w-(r^{-1}s)^{\frac{1}{2}}z_3)\\
&&+(rs)^{\frac{1}{2}}(r{+}(rs)^{\frac{1}{2}}+s)\,(z_1-(r^{-1}s)^{\frac{1}{2}}w)(z_2-(r^{-1}s)^{\frac{1}{2}}w)(w-(r^{-1}s)^{\frac{1}{2}}z_3)\\
&&-(rs)^{\frac{3}{2}}(z_1-(r^{-1}s)^{\frac{1}{2}}w)(z_2-(r^{-1}s)^{\frac{1}{2}}w)(z_3-(r^{-1}s)^{\frac{1}{2}}w) \Big)=0,
\end{eqnarray*}
which is almost the same as that of the one-parameter case with $q=(rs^{-1})^\frac{1}{2}$ (\cite{JKM2}). 
Thus we complete the proof of Theorem 4.1.

\medskip

\vskip30pt \centerline{\bf ACKNOWLEDGMENT}

\bigskip
N. Hu was supported in part by the NNSF of China (No. 11271131).  H. Zhang would
like to thank the support of NNSF of China (No. 11371238).

\bigskip
\def\refname{\cen{\bf REFERENCES}}

\end{document}